\documentclass{amsart}

\usepackage{amsfonts}
\usepackage{amssymb}

\newcommand\R{{\mathbf{R}}}

\newcommand\E{{\mathbf{E}}}
\newcommand\B{{\mathbf{B}}}
\newcommand\Z{{\mathbf{Z}}}
\newcommand\M{{\mathbf{M}}}
\newcommand\I{{\mathcal{I}}}
\newcommand\X{{\mathcal{X}}}
\newcommand\Y{{\mathcal{Y}}}

\newcommand\eps{\varepsilon}
\newcommand\pt{\operatorname{pt}}
\newcommand\N{{\mathbf{N}}}

\theoremstyle{plain}
  \newtheorem{theorem}[subsection]{Theorem}

  \newtheorem{proposition}[subsection]{Proposition}
  
  \newtheorem{lemma}[subsection]{Lemma}

\theoremstyle{remark}
  \newtheorem{remark}[subsection]{Remark}
  
  \newtheorem{example}[subsection]{Example}

\theoremstyle{definition}
  \newtheorem{definition}[subsection]{Definition}

\include{psfig}

\begin{document}

\title[Convergence of multiple ergodic averages]{Norm convergence of multiple ergodic averages for commuting transformations}

\author{Terence Tao}
\address{UCLA Department of Mathematics, Los Angeles, CA 90095-1596.}
\email{tao@@math.ucla.edu}

\begin{abstract}  Let $T_1, \ldots, T_l: X \to X$ be commuting measure-preserving transformations on a probability space $(X, \X, \mu)$.  We show that the multiple ergodic averages $\frac{1}{N} \sum_{n=0}^{N-1} f_1(T_1^n x) \ldots f_l(T_l^n x)$ are convergent in $L^2(X,\X,\mu)$ as $N \to \infty$ for all $f_1,\ldots,f_l \in L^\infty(X,\X,\mu)$; this was previously established for $l=2$ by Conze and Lesigne \cite{CL} and for general $l$ assuming some additional ergodicity hypotheses on the maps $T_i$ and $T_i T_j^{-1}$ by Frantzikinakis and Kra \cite{fkra} (with the $l=3$ case of this result established earlier in \cite{zhang}).  Our approach is combinatorial and finitary in nature, inspired by recent developments regarding the hypergraph regularity and removal lemmas, although we will not need the full strength of those lemmas.  In particular, the $l=2$ case of our arguments are a finitary analogue of those in \cite{CL}.
\end{abstract}

\maketitle

\section{Introduction}

The purpose of this paper is to establish the following norm convergence result for multiple commuting transformations.

\begin{theorem}[Norm convergence]\label{main}  Let $l \geq 1$ be an integer.  Assume that $T_1,\ldots,T_l: X \to X$ are commuting invertible measure-preserving transformations of a measure space $(X, \X, \mu)$.  Then for any $f_1,\ldots,f_l \in L^\infty(X,\X,\mu)$, the averages
$$\frac{1}{N} \sum_{n=0}^{N-1} f_1(T_1^n x) \ldots f_l(T_l^n x)$$
are convergent in $L^2(X,\X,\mu)$.
\end{theorem}

\begin{remark} By using H\"older's inequality and a limiting argument, one can relax the $L^\infty$ conditions on $f_i$ to $L^{p_i}$ conditions for certain finite exponents $p_i$.  For similar reasons, one can also replace the $L^2$ norm with other $L^p$ norms, provided that $\frac{1}{p} \geq \frac{1}{p_1} + \ldots + \frac{1}{p_l}$.  We omit the standard details.
\end{remark}

The case $l=1$ is essentially the mean ergodic theorem.  The case $l=2$ is due to Conze and Lesigne \cite{CL}.  This result had been obtained by Zhang \cite{zhang} for $l=3$ and Frantzikinakis and Kra \cite{fkra} for general $l$ under the additional hypotheses that each of the $T_i$ and the $T_i T_j^{-1}$ (for $i \neq j$) are \emph{individually} ergodic transformations.  The result was also obtained by Lesigne \cite{les} for certain distal systems and by Berend and Bergelson \cite{berend} for certain weakly mixing systems.  In the important special case $T_i = T^i$ for some measure-preserving transformation $T: X \to X$, this result was first obtained for general $l$ by Host and Kra \cite{host-kra2} (with a different proof given subsequently by Ziegler \cite{ziegler}).

All of the preceding arguments mentioned above approach the norm convergence problem through the techniques of ergodic theory, for instance by constructing characteristic factors for the above system.  Here we shall adopt a somewhat different-looking approach, which is based on running the Furstenberg correspondence principle \emph{in reverse} to deduce the above ergodic theory result from a purely combinatorial result (much as the Furstenberg recurrence theorem \cite{furst} can be deduced from Szemer\'edi's theorem \cite{szemeredi}).  More precisely, we shall deduce Theorem \ref{main} from the following ``finitary'' version, in which the general measure-preserving system $(X,\X,\mu,T_1,\ldots,T_l)$ has been replaced by the finite abelian group $\Z_P^l = (\Z/P\Z)^l$ for some large integer $P$, with the discrete $\sigma$-algebra, the uniform probability measure, and the standard $l$ commuting shifts $T_i x := x + e_i$.

\begin{definition}[Expectation notation]\label{expect}  For any finite set $B$ and any function $f: B \to \R$, we use $|B|$ to denote the cardinality of $B$, and define the average $\E_{x \in B} f(x) := \frac{1}{|B|} \sum_{x \in B} f(x)$.  In particular, if $N$ is a positive integer, we use $[N]$ to denote the discrete interval $[N] := \{ 0,1,\ldots,N-1\}$, and thus $\E_{n \in [N]} f(n) = \frac{1}{N} \sum_{n=0}^{N-1} f(n)$.
\end{definition}

\begin{definition}[Finitary averages]\label{finat}  Let $l \geq 1$ and $P \geq 1$.  We let $e_1,\ldots,e_l$ be the standard generators of the finite additive group $\Z_P^l$.  For any functions $f_1, \ldots, f_l: \Z_P^l \to \R$ and any $N \geq 1$, we define the multiple average $A_N(f_1,\ldots,f_l): \Z_P^l \to \R$ by the formula
$$ A_N(f_1,\ldots,f_l)(a) := \E_{n \in [N]} \prod_{i=1}^l f_i(a + n e_i).$$
\end{definition}

\begin{example}\label{anfn}  If $l=2$ and $f_1,f_2: \Z_P^2 \to \R$, then
$$ A_N(f_1,f_2)(v_1,v_2) = \frac{1}{N} \sum_{n=0}^{N-1} f_1(v_1+n,v_2) f_2(v_1,v_2+n)$$
for all $v_1,v_2 \in \Z_P$.
\end{example}

We let $\N := \{1,2,3,\ldots\}$ denote the positive natural numbers.

\begin{theorem}[Finitary norm convergence]\label{main2}  Let $l \geq 1$ be an integer, let $F: \N \to \N$ be a function, and let $\eps > 0$.  Then there exists an integer $M^* > 0$ with the following property:  If $P \geq 1$ and $f_1,\ldots,f_l: \Z_P^l \to [-1,1]$ are arbitrary functions on $\Z_P^l$, then there exists an integer $1 \leq M \leq M^*$ such that we have the ``$L^2$ metastability''
\begin{equation}\label{l2-stability}
 \| A_N(f_1,\ldots,f_l)  - A_{N'}(f_1,\ldots,f_l) \|_{L^2(\Z_P^l)} \leq \eps
\end{equation}
for all $M \leq N, N' \leq F(M)$, where we give $\Z_P^l$ the uniform probability measure.
\end{theorem}

\begin{remark} For applications, Theorem \ref{main2} is only of interest in the regime where $F(M)$ is much larger than $M$, and $P$ is extremely large compared to $l$, $F$, or $\eps$.  The key points are that the function $F$ is arbitrary (thus one has arbitrarily high quality regions of $L^2$ metastability), and that the upper bound $M^*$ on $M$ is independent of $P$.  The $l=1$ version of this theorem was essentially established (with $\Z_P^l$ replaced by an arbitrary measure-preserving system) in \cite{agt}.
\end{remark}

\begin{remark} The presence of the arbitrary function $F: \N \to \N$ may appear strange, but this is in fact a natural consequence of the ``quantifier elimination'' necessary\footnote{In proof theory, this finitisation is known as the \emph{G\"odel functional interpretation} of the infinitary statement, which is also closely related to the \emph{Kriesel no-counterexample interpretation} \cite{kriesel}, \cite{kriesel2} or \emph{Herbrand normal form} of such statements; see \cite{ulrich} for further discussion.  We thank Ulrich Kohlenbach for pointing out this connection.} in order to finitise a convergence result.  For instance, if $f_1,f_2,\ldots$ are a sequence in a normed vector space $V$, observe that the statement
$$ f_1,f_2,\ldots \hbox{ are a Cauchy sequence in } V$$
is by definition equivalent to the assertion that for every $\eps > 0$ there exists $M \geq 1$ such that
$$ \| f_N - f_{N'} \|_V \leq \eps \hbox{ for all } N, N' \geq M,$$
and that this in turn is equivalent to the assertion that for every $\eps > 0$ and every $F: \N \to \N$, there exists $M \geq 1$ such that
\begin{equation}\label{fnv}
 \| f_N - f_{N'} \|_V \leq \eps \hbox{ for all } M \leq N, N' \leq F(M).
\end{equation}
Philosophically, the statement \eqref{fnv} looks easier to prove because (once one fixes the function $F$) one is only asking for the sequence $f_N$ to be \emph{metastable} rather than \emph{stable} - i.e. stable on a finite range $[M,F(M)]$ rather than an infinite range $[M,+\infty)$.  This allows us to perform pigeonholing tricks based on locating several disjoint intervals of the form $[M,F(M)]$, as was recently carried out in \cite{tao:besi}.  Indeed our arguments here have some of the ``multiscale analysis'' flavour of \cite{tao:besi}.  See also \cite{tao:hyper}, in which functions such as $F$ play a key role in establishing a hypergraph regularity lemma.
\end{remark}

We shall establish Theorem \ref{main2} by ``finitary ergodic theory'' techniques, reminiscent of those used in \cite{gt-primes} to establish arbitrarily long arithmetic progressions in the primes.  For instance, instead of building infinitary characteristic factors as was done in earlier work on this problem, we shall build finitary characteristic factors out of ``anti-uniform functions'' analogous to those in \cite{gt-primes}.  This allows us to essentially reduce Theorem \ref{main} to the case in which all the functions $f_1,\ldots,f_l$ are anti-uniform functions (which will in turn be polynomial combinations of basic anti-uniform functions).  The anti-uniformity allows one to reduce the complexity of the average, and very roughly speaking allows one to deduce the $l$-dimensional convergence result in Theorem \ref{main2} from an $l-1$-dimensional convergence result\footnote{This is analogous to how the argument in \cite{CL} deduced the $l=2$ case of Theorem \ref{main} from various one-dimensional convergence results such as the mean and Birkhoff ergodic theorems.  Indeed our own proof of the $l=2$ case of Theorem \ref{main} was inspired (albeit somewhat indirectly) by the arguments in \cite{CL}.}.  However, for technical reasons, we will not induct on Theorem \ref{main2} directly, but on a more complicated counterpart (see Theorem \ref{main3} below), and induct on a ``complexity'' $d$ rather than a ``dimension'' $l$.

Interestingly, the theory of nilsystems (or spectral theory, or Fourier analysis) does not play any role in our arguments (in sharp contrast to \cite{host-kra2} or \cite{ziegler}), although the cubes and Gowers-type norms which appear in \cite{host-kra2} have a faint presence here via our machinery of anti-uniform functions.  Similarly, the full strength of tools such as the hypergraph regularity lemma are not needed; instead we need the weaker ``Koopman-von Neumann'' counterparts to such regularity lemmas (analogous to the ``weak regularity lemma'' of Frieze and Kannan \cite{fkannan}).  As with other applications of graph and hypergraph methods, the $\Z^l$ group action in fact plays remarkably little role in these arguments, although the standard fact that this group is amenable\footnote{For instance, one can establish analogues of our results in which $\Z$ is replaced with the infinite vector space $F^\N$ over a finite field $F$ generated by an infinite basis $e_1,e_2,\ldots$, and the intervals $[N]$ are replaced with the subspaces spanned by $e_1,\ldots,e_N$.  In fact the proof in this finite field case is somewhat easier than in the integer case due to the perfect nesting of the scales.} will be implicitly used at several crucial junctures (basically allowing us to treat coarse scales averages as an average of fine scale averages, modulo negligible errors\footnote{For a specific example of this, if $T$ is a shift operator and $S_N$ are the averaging operators $S_N := \E_{n \in [N]} T^n$, observe for $1 \leq M \leq N$ that $S_M S_N$ and $S_N$ differ in $L^2$ operator norm by only $O(M/N)$, and thus we have the heuristic $S_M S_N \approx S_N$ in the regime $N \gg M$.}).

The main advantage of working in the finitary setting, as opposed to the more traditional infinitary one, is that the underlying dynamical system becomes extremely explicit, being simply the standard shifts on $\Z_P^l$.  In particular we have a Cartesian product structure which allows us to construct various product sets\footnote{Actually, as is usual in the hypergraph approach to recurrence problems, we shall lift $\Z_P^l$ to $\Z_P^{l+1}$ in order to abstract away the arithmetic aspects of the shift operations; see Section \ref{gensec} below.} in our dynamical system, without having to pay attention to technical issues such as measurability.  This product structure will be crucial to our arguments (it basically endows our system with the structure of a hypergraph).  It seems of interest to try to obtain similar product structures in the traditional infinitary setting; the argument in \cite{CL} achieves this to some extent in the $l=2$ case. (See also \cite{tao:corr} for another (not entirely satisfactory) attempt to endow dynamical systems with hypergraph structure.)  This would likely lead to a more traditional infinitary proof of Theorem \ref{main}.  

Finally, we remark that our methods do not seem to extend to the significantly more difficult question of pointwise almost everywhere convergence of these averages; for that task, some sort of multilinear maximal inequality may be needed.

\subsection{Organisation}

This paper is organised as follows.  Firstly, in Section \ref{revsec}, we use the Furstenberg correspondence principle in the reverse direction to deduce the infinitary convergence theorem, Theorem \ref{main}, from its finitary counterpart, Theorem \ref{main2}.  Then, in Section \ref{notation-sec}, we set out our basic notation we need to establish Theorem \ref{main2}.  In Section \ref{gensec}, we deduce Theorem \ref{main2} from a more technical variant, Theorem \ref{main3}, which is in a form suitable for applying mathematical induction on a certain ``complexity'' parameter $d$.  The base case $d=1$ (which is essentially a finitary analogue of the mean ergodic theorem, as in \cite{agt}) is then handled in Section \ref{base-sec}; these arguments are then generalised to handle the inductive case $d>1$ in Section \ref{inductive-sec}.  

At several points in the argument it will be convenient to pass from a ``probabilistic'' norm convergence result to a ``deterministic'' one.  The natural tool for this is the Lebesgue dominated convergence theorem, but as we are working in a finitary setting, we of course need a finitary counterpart of this infinitary theorem.  We discuss such a counterpart in an appendix to this paper.  Actually, it is possible to skip this dominated convergence step and work entirely in a probabilistic setting throughout, but this causes the notation to be slightly more complicated.

\subsection{Acknowledgements}

We thank Ciprian Demeter for explaining the argument in \cite{CL}, for encouragement, and for bringing the norm convergence problem to our attention.  We thank Henry Towsner and Ulrich Kohlenbach for bringing the author's attention to \cite{agt} and to pointing out the connections to proof theory.  We also thank Jennifer Chayes for suggesting the term ``metastability'', and Tim Austin, Ciprian Demeter, Henry Towsner, and Christoph Thiele for helpful discussions.  Finally, we thank the anonymous referee for a careful reading of the manuscript and for many suggestions and corrections.  The author is supported by NSF grant CCF-0649473 and a grant from the MacArthur Foundation.

\section{The reverse Furstenberg correspondence principle}\label{revsec}

In this section we show how to reverse the Furstenberg correspondence principle \cite{furst} to deduce Theorem \ref{main} from Theorem \ref{main2}.

\begin{proof}[Proof of Theorem \ref{main} assuming Theorem \ref{main2}]
Observe that the $l$ commuting transformations generate a measure-preserving action of $\Z^l$ on the system $(X,\X,\mu)$.  We claim that we may reduce\footnote{Actually, this reduction step, as well as the step involving the generic point $x_0$ below, is not strictly necessary to our argument, provided that one is willing to replace Theorem \ref{main2} by the more complicated-looking generalisation in Theorem \ref{main3} below.} to the case when this action is ergodic, i.e. the only sets which are invariant under \emph{all} of the $T_1,\ldots,T_l$ have either zero measure or full measure.  Note that this is a much weaker property than requiring that each of the $T_1,\ldots,T_l$ (or the $T_i T_j^{-1}$) are \emph{individually} ergodic.  This reduction is standard and performed for instance in \cite[page 157]{CL}, so we only sketch it here.  Using the ergodic decomposition (see e.g. \cite{furst-book}) one can disintegrate $\mu$ as an integral of measures $\mu_y$, such that each $\mu_y$ is invariant and ergodic with respect to the $\Z^l$ action.  By hypothesis, the averages
$\E_{n \in [N]} f_1(T_1^n x) \ldots f_l(T_l^n x)$ are convergent, hence Cauchy, in each of the $L^2(X,\X,\mu_y)$; they are also bounded between $-1$ and $1$.  By the dominated convergence theorem we conclude that these averages are Cauchy, hence convergent, in $L^2(X,\X,\mu)$, as desired.

Henceforth we assume the $\Z^l$ action to be ergodic on $(X,\X,\mu)$.  Applying the Birkhoff pointwise ergodic theorem for $\Z^l$ (see \cite{wiener}), and in particular we see that for any $f \in L^\infty(X,\X,\mu)$ that
\begin{equation}\label{pox} \lim_{P \to \infty} \E_{v \in [P]^l} f( T^v x_0 ) = \int_X f\ d\mu
\end{equation}
for almost every $x_0$, where we adopt the convention
$$ T^{(v_1,\ldots,v_l)} := T_1^{v_1} \ldots T_l^{v_l}.$$
Let us say that a point $x_0$ is \emph{generic} if \eqref{pox} holds for all functions $f$ which are polynomial combinations of the $f_1,\ldots,f_l$ with rational coefficients.  Since there are only countably many such functions, we see that almost every point is generic.

Fix a generic point $x_0$.  Recall that our objective is to show that the sequence of functions
$$\E_{n \in [N]} \prod_{i=1}^l f_i(T_i^n x)$$
is convergent in $L^2(X,\X,\mu)$.  It of course suffices to show that it is a Cauchy sequence.  If this is not the case, then there exists $\eps > 0$ with the property that for every integer $M > 0$, there exists an integer $F(M) > M$ such that
\begin{equation}\label{xfm}
 \int_X \left| \E_{n \in [F(M)]} \prod_{i=1}^l f_i(T_i^n x)
- \E_{n \in [M]} \prod_{i=1}^l f_i(T_i^n x)\right|^2\ d\mu(x) > 3\eps^2
\end{equation}
(say).  Fix this $\eps$ and $F$.  Applying \eqref{pox}, we can write the left-hand side of \eqref{xfm} as
$$ \lim_{P \to \infty} \E_{v \in [P]^l}
\left|\E_{n \in [F(M)]} \prod_{i=1}^l f_i(T_i^n T^v x_0)
- \E_{n \in [M]} \prod_{i=1}^l f_i(T_1^n T^v x_0)\right|^2.$$
Let $M^*$ be the integer depending on $l, \eps, F$ which appears in Theorem \ref{main2}.  
Then, if $P$ is sufficiently large depending on $M^*, F, f_1, \ldots, f_l, x_0, \eps$, we can ensure that
\begin{equation}\label{snapple}
\E_{v \in [P]^l}
\left|\E_{n \in [F(M)]} \prod_{i=1}^l f_i(T_i^n T^v x_0)
- \E_{n \in [M]} \prod_{i=1}^l f_i(T_i^n T^v x_0)\right|^2 > 2\eps^2
\end{equation}
for all $1 \leq M \leq M^*$.

We now assume $P$ large enough so that the above properties hold.  Define the functions $g_1,\ldots,g_l: \Z_P^l \to [-1,1]$ by setting
$$ g_i(v) := f_i( T^v x_0 )$$
for all $v \in \Z_P^l$, where we artificially identify $\Z_P$ with $[P]$ in the usual manner.  From \eqref{snapple} we see that
$$ \E_{v \in [P]^l} \left|A_{F(M)}(g_1,\ldots,g_l)(v) - A_{M}(g_1,\ldots,g_l)(v)\right|^2 > \eps^2$$
for all $1 \leq M \leq M^*$, if $P$ is large enough depending on $M^*, F, \eps$ (this is necessary to be able to neglect the (rare) ``wraparound effects'' caused when the shifts $T_1^n,\ldots,T_l^n$ push one of the coefficients of $a$ beyond $P$).  But this contradicts Theorem \ref{main2}.  This contradiction establishes Theorem \ref{main} as desired.
\end{proof}

\begin{remark}  It is also possible to apply the Furstenberg correspondence principle (as in \cite{furst} or \cite{furst-book}) in the more standard direction and deduce Theorem \ref{main2} from Theorem \ref{main}, by using the weak sequential compactness of probability measures on the discrete cube $\{0,1\}^{\Z^l}$.  We leave the details to the interested reader.
\end{remark}

It remains to prove Theorem \ref{main2}.  This will be the purpose of later sections.

\section{Finitary notation}\label{notation-sec}

Theorem \ref{main2} is a statement in ``finitary'' mathematics - it concerns averages over finite sets.  In this section we lay out some finitary notation which will be of use in establishing that theorem (and also point out some connections with graph and hypergraph theory which are implicitly lurking just beneath the surface).  We will of course be heavily using the expectation notation in Definition \ref{expect}. We also recall some standard asymptotic notation:

\begin{definition}[Asymptotic notation] We use $A \ll B$ or $B \gg A$ to denote the bound $A \leq CB$ for some constant $C$, and $O(A)$ to denote any quantity bounded in magnitude by $CA$.  If we wish to allow the constant $C$ to depend on auxiliary parameters, we will denote this by subscripts, e.g. $O_\eta(A)$ denotes a quantity bounded by $C_\eta A$ where $C_\eta$ is allowed to depend on $\eta$.
\end{definition}

\subsection{Factors}

Next, we recall the notion of a factor from ergodic theory.

\begin{definition}[Factor]  Let $(X,\X,\mu)$ be a probability space.  A \emph{factor} of $(X,\X,\mu)$ is a triplet $\Y = (Y,\Y,\pi)$, where $Y$ is a set, $\Y$ is a $\sigma$-algebra, and $\pi: X \to Y$ is a measurable map.  If $\Y$ is a factor, we let $\B_\Y := \{ \pi^{-1}(E): E \in \Y \}$ be the sub-$\sigma$-algebra of $\X$ formed by pulling back $\Y$ by $\pi$.  A function $f: X \to \R$ is said to be \emph{$\Y$-measurable} if it is measurable with respect to $\B_\Y$.  If $f \in L^2(X,\X,\mu)$, we let $\E(f|Y) = \E(f|\B_Y)$ be the orthogonal projection of $f$ to the closed subspace $L^2(X,\B_Y,\mu)$ of $L^2(X,\X,\mu)$ consisting of $\Y$-measurable functions.  If $\Y = (Y,\Y,\pi)$ and $\Y' = (Y',\Y',\pi')$ are two factors, we let $\Y \vee \Y'$ denote the factor $(Y \times Y', \Y \otimes \Y', \pi \oplus \pi')$.
\end{definition}

\begin{remark} The concept of a factor in ergodic theory corresponds closely with the concept of a \emph{partition} or \emph{colouring} in graph or hypergraph theory.
\end{remark}

Our probability spaces shall usually be finite sets with the uniform distribution.  More precisely,
if $Y$ is a finite set, let $2^Y = \{ E: E \subset Y \}$ be the discrete $\sigma$-algebra on $Y$, and let $\mu_Y$ be the uniform probability measure on $Y$.  In particular,  $L^2(Y)$ be the finite-dimensional real Hilbert space of functions $f: Y \to \R$, endowed with the inner product
$$ \langle f, g \rangle_{L^2(Y)} := \E_{y \in Y} f(y) g(y).$$
More generally, if $X = (X,\X,\mu)$ is another probability space, $L^2(Y \times X)$ is the real Hilbert space of measurable functions $f: Y \times X \to \R$, endowed with the inner product
$$ \langle f, g \rangle_{L^2(Y \times X)} := \int_X \E_{y \in Y} f(y,x) g(y,x)\ d\mu(x).$$

\begin{remark} Our use of the uniform distribution for $Y$ corresponds to the customary convention in graph and hypergraph theory to give all vertices, edges, etc. equal weight.  One can of course replace uniform distributions by more general probability distributions, corresponding to weighted graphs and hypergraphs, but we will not need to do so here.  
\end{remark}

In the infinitary theory, we can use any measurable function $f: X \to \R$ to generate a factor of $X$, whose $\sigma$-algebra is generated by the level sets $f^{-1}([a,b])$ for all $a,b$.  The function $f$ will then be measurable with respect to that factor.  Such factors turn out to be far too large for us to use in the finitary setting (for instance, if $X$ is finite and $f$ takes different values at each point of $X$, then the above $\sigma$-algebra will be the maximal $\sigma$-algebra $2^X$).  Instead, we will need some slightly coarser factors, defined via the following lemma.

\begin{lemma}[Each function generates its own factor]\label{factor-lemma}  Let $(X,\X,\mu)$ be a probability space, let $I \subset \R$ be a compact interval, and let $\varphi: X \to I$ be a measurable function.  Then for any $\eta_0 > 0$ there exists a factor $\Y_{\eta_0}(\varphi)$ with the following properties.
\begin{itemize}
\item[(i)] ($\varphi$ lies in its own factor) For any factor $\Y'$, we have
$$ \| \varphi - \E(\varphi|\Y_{\eta_0}(\varphi) \vee \Y') \|_{L^\infty(X,\X,\mu)} \leq \eta_0.$$
\item[(ii)] (Bounded complexity)  The $\sigma$-algebra $\B_\Y$ is generated by $O_{I,\eta_0}(1)$ atoms.
\item[(iii)] (Approximation by polynomials of $\varphi$) If $A$ is any atom in $\B_\Y$ and $\eta_1 > 0$, there exists a polynomial $\Psi_A: I \to [0,1]$ of degree $O_{I,\eta_0,\eta_1}(1)$ and coefficients $O_{I,\eta_0,\eta_1}(1)$ such that
$$ \| 1_A - \Psi_A(\varphi) \|_{L^1(X,\X,\mu)} \leq \eta_1$$
and
$$ \| 1_A - \Psi_A(\varphi) \|_{L^\infty(X,\X,\mu)} \leq 1.$$
\end{itemize}
\end{lemma}

\begin{proof}  This lemma essentially appears in \cite[Proposition 7.2]{gt-primes}, \cite[Proposition 6.1]{tao:gauss}, or \cite[Proposition 7.3]{tz}, so we give only a brief sketch here.

We use the probabilistic method.  Let $\alpha \in [0,1]$ be chosen uniformly at random.  We let $\Y(\varphi) = \Y_\alpha(\varphi)$ be the factor
$$ \Y(\varphi) = (I, \B_{\alpha,\eta_0}, \varphi)$$
where $\B_{\alpha,\eta_0}$ is the $\sigma$-algebra generated by the intervals $[(n+\alpha+1)\eta_0, (n+\alpha)\eta_0)$.  The properties (i), (ii) are then obvious, so it suffices to verify (iii).  Firstly, we observe that it suffices to verify (iii) in the case where $\eta_1 = 2^{-j}$ for an integer $j \geq 0$.  We will in fact show that for each fixed $j$, that (iii) holds with probability $1 - O_{I,\eta_0}(2^{-j})$; from the union bound we thus see that there exists a choice of $\alpha$ for which (iii) holds for all $j$ that are sufficiently large depending on $I,\eta_0$, and the claim then follows since the claim for small $j$ clearly follows from that of large $j$.

Let us now fix $j$.  By (ii) and the union bound again, it suffices to verify the claim for a single atom $A = \varphi^{-1}([(n+\alpha+1)\eta_0, (n+\alpha)\eta_0))$, where $n \in \Z$ is fixed.  We define the exceptional set
$$ B := \{ x \in X: |\varphi(x) - (n-\alpha-i)\eta_0| \leq 2^{-2j} \hbox{ for some } i=0,1 \}$$
then from Fubini's theorem we see that $B$ has small measure on the average:
$$ \E \mu(B) \ll 2^{-2j}.$$
By Markov's inequality, we thus see that $\mu(B) \leq 2^{-j}/2$ with probability $1 - O(2^{-j})$.  We now apply Urysohn's lemma followed by the Weierstrass approximation theorem to locate a polynomial $\Psi_A: I \to [0,1]$ of degree $O_{I,\eta_0,j}(1)$ and coefficients $O_{I,\eta_0,j}(1)$ such that
$$ |\Psi_A(t) - 1_{[(n+\alpha+1)\eta_0, (n+\alpha)\eta_0)}(t)| \leq 2^{-j}/2$$
for all $t$ with $|t - (n-\alpha-i)\eta_0| \geq 2^{-2j}$ for $i=0,1$.  (Note that $\alpha$ ranges in a compact set, and so the bounds on the degree and coefficients on $\Psi_A$ are uniform in $\alpha$.)  One then easily verifies that
$$ \| 1_A - \Psi_A(\varphi) \|_{L^1(X,\X,\mu)} \leq 2^{-j}/2 + \mu(B) \leq 2^{-j}$$
and the claim (iii) follows.
\end{proof}

Henceforth we fix the assignment $(\varphi, \eta_0) \mapsto \Y_{\eta_0}(\varphi)$ of a factor to each function $\varphi$ and an error tolerance $\eta_0$ as defined above.

\subsection{Products, edge factors, and complexity}

We shall work frequently with finite Cartesian products
\begin{equation}\label{yi}
 Y_I = \prod_{i \in I} Y_i := \left\{ (y_i)_{i \in I}: y_i \hbox{ for all } i \in I \right\}
\end{equation}
where $I$ is a finite index set, and the $Y_i$ are also finite sets.  We of course adopt the usual convention that
$$ Y^n := \prod_{i \in \{1,\ldots,n\}} Y$$
for any non-negative integer $n$.

For technical reasons (basically due to our use of probabilistic methods) we will also need to deal with the slightly larger product spaces
\begin{equation}\label{yix}
Y_I \times X = \prod_{i \in I} Y_i \times X
\end{equation}
where $X = (X,\X,\mu)$ is another probability space (possibly infinite).  The space $X$ should be thought of as a ``passive'' space, as the parameters in $X$ will simply be averaged over at the end of the day, with no non-trivial interaction with any other parameters in the argument.  The space $Y_I \times X$ is then also a probablity space, with the product $\sigma$-algebra $2^{Y_I} \otimes \X$ and the product measure $\mu_{Y_I} \times \mu$.  Of course one can view ordinary Cartesian products \eqref{yi} as a special case of \eqref{yix} in which the probability space $X$ is just a point, $X = \pt$.

\begin{remark} In the graph and hypergraph theory language, the sets $Y_i$ should be viewed as disjoint classes of vertices, and various subsets of $Y_I$ should be interpreted as partite graphs or hypergraphs, where the edges consist of up to one vertex from each of the classes $Y_i$.  Subsets of the larger space $Y_I \times X$ should be interpreted as \emph{random} partite graphs or hypergraphs.
\end{remark}

Now we come to a crucial concept in our product space analysis.

\begin{definition}[Edge factors]\label{edge}  Let $Y_I \times X = (Y_I \times X, 2^{Y_I} \otimes \X, \mu_{Y_I} \times \mu)$ be a probability space as above.  For any $e \subset I$, let $Y_e := \prod_{i \in e} Y_i$, and let $\pi_e: Y_I \times X \to Y_e \times X$ be the edge projection
$$ \pi_e( (y_i)_{i \in I}, x ) := ( (y_i)_{i \in e}, x ).$$
We then let $\Y_e$ be the factor $(Y_e \times X, 2^{Y_e} \times \X, \mu_{Y_e} \times \mu)$ of $Y_I \times \X$.  We say that a function $f: Y_I \times X \to \R$ is \emph{$e$-measurable} if it is $\Y_e$-measurable.
\end{definition}

\begin{example}\label{finex}  Let $Y$ be a finite set, let $X = (X,\X,\mu)$ be a probability space, and let $f: Y^3 \times X \to \R$ be a measurable function.  Then $f$ is $\{1,3\}$-measurable if and only if it takes the form
$$ f(y_1,y_2,y_3,x) = f_{13}(y_1,y_3,x)$$
for some measurable $f_{13}: Y^{\{1,3\}} \times X \to \R$.  Similarly, $f$ is $\{3\}$-measurable if and only if it takes the form
$$ f(y_1,y_2,y_3,x) = f_{3}(y_3,x)$$
for some measurable $f_3: Y^{\{3\}} \times X \to \R$.
\end{example}

\begin{remark} In the graph and hypergraph theory language, an $e$-measurable set should be regarded as an $|e|$-uniform, $e$-partite hypergraph on the vertex classes $Y_i$ for $i \in e$.  For instance, continuing the above example, we let $Y_1, Y_2, Y_3$ be three identical copies of $Y$, viewed as vertex sets, then if an indicator function $1_{E_{13}}$ is $\{1,3\}$-measurable then it can be viewed as describing a bipartite graph connecting $Y_1$ and $Y_3$, whereas if an indicator function $1_{E_3}$ is measurable it can be viewed as describing a set of vertices in $Y_3$.  Finally, a $\{1,2,3\}$-measurable indicator $1_{E_{123}}$ can be viewed as a $3$-uniform tripartite hypergraph 
connecting $Y_1$, $Y_2$, and $Y_3$.
\end{remark}

We make the trivial remark that an $e$-measurable function is automatically $e'$-measurable for any $e' \supset e$.  For instance, all functions are $I$-measurable.

Let $d \geq 1$ be an integer.  We will informally refer to an edge factor $\Y_e$ as having \emph{complexity} $d$ if $|e|=d$.  We would like to combine together all the edge factors $\Y_e$ of a given complexity $d$ to create a ``complexity $d$ factor'', which should morally form a tower of factors in $d$ analogous to the Furstenberg tower constructed for instance in \cite{furst-book}.  However, one has to take some care with this, since as $\sigma$-algebras (or even as algebras), the edge factors $\Y_e$ of complexity $d$ generate the entire $\sigma$-algebra $2^{Y_I} \otimes \X$.  To obtain a meaningful concept of a ``complexity $d$ factor'', then, we have to also limit the complexity of the polynomial combinations of $e$-complexity functions we shall employ.  This leads to the following important definitions.

\begin{definition}[Complexity]\label{complex-def}  Let $Y_I \times X = (Y_I \times X, 2^{Y_I} \otimes \X, \mu_{Y_I} \times \mu)$ be a probability space as above.  Let $1 \leq d \leq |I|$.  A function $g: \prod_{i \in I} Y_i \times X \to [-1,1]$ is a \emph{primitive function of complexity at most $d$} if it takes values in $[-1,1]$ and is $e$-measurable for some $e \subset I$ with $|e| \leq d$.  A function $g: \prod_{i \in I} Y_i \times X \to [-1,1]$ is a \emph{basic function of complexity at most $d$} if it is the product of finitely many primitive functions of complexity at most $d$, or equivalently if it has a representation $g = \prod_{e \subset I: |e|=d} g_e$ where each $g_e$ is $e$-measurable.  A function $g: \prod_{i \in I} Y_i \times X \to \R$ is an \emph{elementary function of complexity at most $(d,J)$} for some integer $J \geq 1$ if it can be expressed as the sum of at most $J$ basic functions of complexity at most $d$.  
\end{definition}

\begin{example} We continue Example \ref{finex}.  If $f_{12}, f_{13}, g_{12}, g_{23}: Y^2 \times X \to [-1,1]$ are measurable functions, then the function
$$ f(y_1,y_2,y_3,x) := f_{12}(y_1,y_2,x)$$
is a primitive function of complexity at most $2$,
$$ f'(y_1,y_2,y_3,x) := f_{12}(y_1,y_2,x) f_{13}(y_1,y_3,x)$$
is a basic function of complexity at most $2$, and
$$ f''(y_1,y_2,y_3,x) := f_{12}(y_1,y_2,x) f_{13}(y_1,y_3,x) + g_{12}(y_1,y_2,x) g_{23}(y_2,y_3,x)$$
is an elementary function of complexity at most $(2,2)$.
\end{example}

\begin{remark}  Observe that if $g$ and $g'$ are elementary functions of complexities at most $(d,J)$ and $(d,J')$ respectively, then $g \pm g'$ and $gg'$ have complexities at most $(d,J+J')$ and $(d,JJ')$ respectively; also, if $\alpha$ is any real number with $|\alpha| \leq L$ for some integer $L$, then $\alpha g$ has complexity at most $(d,JL)$.  Thus the space of functions of bounded complexity is morally an algebra.
\end{remark}

\subsection{Group structure}

Graph and hypergraph theory takes place on vertex sets $Y$ which have no algebraic structure.  However, in our application these sets arise from $\Z$ and will have two additional structures: the additive group structure, and the F{\o}lner-type structure coming from the sets $[N]$ that one is averaging over.  To handle these structures we introduce two useful notations.

\begin{definition}[Summation]  Let $G = (G,+)$ be an additive group, and let $G^I$ be any finite Cartesian power of $G$.  Given any vector $v = (v_i)_{i \in I} \in G^I$, we define the sum $\Sigma(v) \in G$ of $v$ by $\Sigma(v) := \sum_{i \in I} v_i$.
\end{definition}

Clearly, $\Sigma$ is a homomorphism from $G^I$ to $G$.  We shall usually apply this notation with $G=\Z_P$ equal to a cyclic group.

To motivate our next definition, we recall the setup in Example \ref{anfn}.  We rewrite $A_N(f_1,f_2)(v_1,v_2)$ as
$$ A_N(f_1,f_2)(v_1,v_2) = \E_{n \in [N]} f_1(-v_2 - (-v_1-v_2-n),v_2) f_2(v_1, -v_1 - (-v_1-v_2-n)).$$
The point of doing this is that we now see that the $f_1$ factor depends only on $v_2$ and $-v_1-v_2-n$, while the $f_2$ factor depends only on $v_1$ and $-v_1-v_2-n$.  To make these dependencies even clearer, we introduce the $\{2,3\}$-measurable function
$$ g_{\{2,3\}}(v_1,v_2,v_3) := f_1( -v_2-v_3, v_2 ) $$
and the $\{1,3\}$-measurable function
$$ g_{\{1,3\}}(v_1,v_2,v_3) := f_2( v_1, -v_1-v_3 )$$
and observe the identity
\begin{equation}\label{adelta-2}
A_N(f_1,f_2)(v_1,v_2) = \E_{n \in [N]} g_{\{2,3\}} g_{\{1,3\}}( v_1, v_2, -v_1-v_2-n ).
\end{equation}
Thus $A_N(f_1,f_2)$ can be viewed as an average of the product of the $\{2,3\}$-measurable function $g_{\{2,3\}}$ and the $\{1,3\}$-measurable function $g_{\{1,3\}}$ along the diagonal region $\{ (v_1,v_2,v_3): v_3 \in -v_1-v_2-[N] \}$.  

More generally, we can represent averages such as $A_N$ as diagonally averaged projections by introducing the following operator.

\begin{definition}[Diagonally averaged projection]\label{dap}  Let $l \geq 1$ and $P \geq 1$ be integers.  Let $(X,\X,\mu)$ be a probability space.  If $f: \Z_P^{l+1} \times X \to \R$ is a measurable function and $N \geq 1$ is an integer, we define the diagonally averaged projection $\Delta_N f: \Z_P^l \times X \to \R$ to be the function
$$ \Delta_N f( v, x) := \E_{n \in [N]} f((v, -\Sigma(v)-n), x )$$
for all $v \in \Z_P^l$ and $x \in X$.
\end{definition}

\begin{remark} The space $X$ is necessary to our argument for technical inductive reasons but should be neglected at a first reading.
\end{remark}

The projection $\Delta_N$ is related to the averages $A_N$ in Definition \ref{finat} by the easily verified identity
\begin{equation}\label{adelta}
A_N(f_1,\ldots,f_l) = \Delta_N( \prod_{i=1}^l g_{\{1,\ldots,l+1\} \backslash \{i\}} )
\end{equation}
for any $f_1,\ldots,f_l: \Z_P^l \to \R$, where for each $1 \leq i \leq l$, the function $g_{\{1,\ldots,l+1\} \backslash \{i\}}: \Z_P^{l+1} \to \R$ is the $\{1,\ldots,l+1\} \backslash \{i\}$-measurable function
$$g_{\{1,\ldots,l+1\} \backslash \{i\}}(v_1,\ldots,v_{l+1}) := f_i( v_1,\ldots,v_{i-1}, -\sum_{1 \leq j \leq l+1: j \neq i} v_j, v_{i+1},\ldots,v_l).$$

One can verify that when $l=2$, that \eqref{adelta} collapses to \eqref{adelta-2}.

\begin{remark} The above elementary arithmetic manipulations are essentially the same manipulations used in the hypergraph approach (see \cite{soly-roth}, \cite{rodl2}, \cite{gowers-reg}, \cite{tao:hyper}) to Szemer\'edi's theorem \cite{szemeredi} or the Furstenberg-Katznelson theorem \cite{fk}, in order to rewrite the problem in a ``hypergraph'' form, by which we mean that the problem now concerns the averages of products of multiple functions, each of which depends on a different set of variables.  (This corresponds to the problem in hypergraph theory of counting the number of instances of a small fixed hypergraph inside a much larger hypergraph.)
\end{remark}

The operator $\Delta_N$ is clearly linear.  For future reference we also observe the module identity
\begin{equation}\label{gih}
 \Delta_N( g_{\{1,\ldots,l\}} h ) = g_{\{1,\ldots,l\}} \Delta_N(h)
\end{equation}
for any $\{1,\ldots,l\}$-measurable $g_{\{1,\ldots,l\}}: \Z_P^{l+1} \times X \to \R$ and any
$h: \Z_P^{l+1} \times X \to \R$, where by abuse of notation we also view the $\{1,\ldots,l\}$-measurable function $g_{\{1,\ldots,l\}}$ as a function on $\Z_P^l$.

\section{A generalisation of Theorem \ref{main}}\label{gensec}

We will prove Theorem \ref{main2} by an induction on the ``complexity'' of the functions $f$ involved.  As it turns out, a naive induction based on Theorem \ref{main2} in its current form does not seem to work well, and so we shall instead use the following more complicated generalisation of Theorem \ref{main2} to induct upon, in which functions such as $f_1,\ldots,f_l$ are allowed to be ``random'' rather than ``deterministic'' (or more precisely, they are allowed to depend on an additional probability space $(X,\X,\mu)$), and have varying levels of ``complexity''.

Specifically, we shall deduce Theorem \ref{main2} from the following more technical variant.

\begin{theorem}[Finitary norm convergence, technical generalisation]\label{main3}  Let $1 \leq d \leq l$, $M_* \geq 1$, and $J \geq 1$ be integers.  Let $F: \N \to \N$ be a function, and let $\eps > 0$.  Then there exists an integer $M^* \geq M_*$ with the following property:  If $P \geq 1$, if $(X,\X,\mu)$ is a probability space, and $g: \Z_P^{l+1} \times X \to \R$ is an elementary function of complexity at most $(d,J)$, then there exists an integer $M_* \leq M \leq M^*$ such that
\begin{equation}\label{l2-stability-tech}
\| \Delta_N(g) - \Delta_{N'}(g) \|_{L^2(\Z_P^l \times X)} \leq \eps
\end{equation}
for all $M \leq N, N' \leq F(M)$.  
\end{theorem}

\begin{remark} This theorem is faintly reminiscent of the ``hypergraph counting lemmas'' which appear for instance in \cite{nrs}, \cite{gowers-reg}, \cite{tao:hyper}.
\end{remark}

The deduction of Theorem \ref{main3} from Theorem \ref{main2} is immediate by specialising to the case where $d=l$ and $M_*=J=1$, where $X$ is a point, and $g$ is the function $\prod_{i=1}^l g_{\{1,\ldots,l+1\} \backslash \{i\}}$, which is a basic function of complexity $d$, and then using \eqref{adelta}.  

\begin{remark} The main point of generalising Theorem \ref{main2} to Theorem \ref{main3} is that it introduces a new parameter $d$ - the maximum \emph{complexity} of all the functions $g_e$ involved.  We shall in fact prove Theorem \ref{main3} by an induction on this parameter $d$ (keeping the dimension $l$ fixed).  The addition of the probability space $(X,\X,\mu)$ is a technical convenience for us, as it allows us to perform a number of averaging or probabilistic arguments without losing the ability to exploit the induction hypothesis.  The passage from one level of complexity $d$ to the next $d+1$ is roughly analogous to that of passing from one dynamical system to a primitive extension in ergodic theory.
\end{remark}

It remains to prove Theorem \ref{main3}.  This will be the purpose of the later sections.

\section{The base case}\label{base-sec}

In this section we shall establish the base case\footnote{In fact, one could incorporate this case into the inductive case, by making $d=0$ the base case instead, but we have chosen to do the $d=1$ case in detail for didactic reasons, as it serves to motivate the higher $d$ argument.} $d=1$ of Theorem \ref{main3}.  

We first make some simple reductions.  Firstly, we can reduce to the case $M_*=1$, by replacing $F(M)$ by the function $\tilde F(M) := F(\max(M,M_*))$, applying Theorem \ref{main3} with $\tilde F$ (and $M_*$ replaced by $1$), and then replacing $M$ with $\max(M,M_*)$.

Next, we reduce to the case $J=1$ by the following argument.  Since $g: \Z_P^{l+1} \times X \to \R$ has complexity at most $(d,J)$, we can write $g = g_1 + \ldots + g_J$ where each $g_k: \Z_P^{l+1} \times X \to \R$ is a basic function of complexity at 
most $d$.  We then define the extended probability space $\tilde X := X \times \{1,\ldots,J\}$, where we give $\{1,\ldots,J\}$ the discrete $\sigma$-algebra and uniform probability measure, and give $\tilde X$ the associated product measure.  We also define the function $\tilde g: \Z_P^{l+1} \times \tilde X \to [-1,1]$ by $\tilde g( v, (x,k) ) := g_k(v,x)$.  One easily verifies from Definition \ref{complex-def} that $\tilde g$ is a basic function of complexity at most $d$, and that we have the identity
$$ \|  \Delta_N(g) - \Delta_{N'}(g) \|_{L^2(\Z_P^l \times X)} 
= J^{1/2} \| \Delta_N(\tilde g) - \Delta_{N'}(\tilde g)\|_{L^2(\Z_P^l \times \tilde X)}
$$
for all $N, N'$.  Because of this, we see that we can reduce to the $J=1$ case (after adjusting $\eps$ by a factor of $J^{1/2}$).

Since $J=1$ and $d=1$, we can now write $g = \prod_{i=1}^{l+1} g_{\{i\}}$ where each $g_{\{i\}}: \Z_P^{l+1} \times X \to \R$ is $\{i\}$-measurable and takes values in $[-1,1]$.  The contributions of the factors $g_{\{i\}}$ with $1 \leq i \leq l$ can be quickly discarded by using the module identity \eqref{gih}. Because of this, we may assume without loss of generality that $\I$ consists only of the singleton set $\{l+1\}$, thus we now just have a single function $g_{\{l+1\}}: \Z_P^{l+1} \times X \to [0,1]$.  We can use the $\{l+1\}$-measurability to write
$$ g_{\{l+1\}}(v_1,\ldots,v_{l+1},x) = g(-v_{l+1},x)$$
where $g: \Z_P \times X \to [-1,1]$ is a measurable function.  We now observe from Definition \ref{dap} that the function $\Delta_N(g_{\{l+1\}})((v_1,\ldots,v_l),x)$ only depends on $v_1+\ldots+v_l$ and $x$.  Thus we may quotient out by the hyperplane $\{ (v_1,\ldots,v_l,v_{l+1}) \in \Z_P^{l+1}: v_1 + \ldots + v_l = 0 \}$ and reduce $\Z_P^{l+1}$ to a one-dimensional group $\Z_P$. We are now reduced to showing the following:

\begin{theorem}[Finitary norm convergence, base case]\label{main3-base}  
Let $F: \N \to \N$ be a function, and let $\eps > 0$.  Then there exists an integer $M^* \ge 1$ with the following property:  If $P \geq 1$, if $(X,\X,\mu)$ is a probability space, and $g: \Z_P \times X \to [0,1]$ is a measurable function, then there exists an integer $1 \leq M \leq M^*$ such that
\begin{equation}\label{l2-stability-base}
\| S_N g - S_{N'} g \|_{L^2(\Z_P \times X)} \leq \eps
\end{equation}
for all $M \leq N, N' \leq F(M)$, where $S_N$ is the averaging operator $S_N g(v,x) := \E_{n \in [N]} g(v+n,x)$, and similarly for $S_{N'}$.  
\end{theorem}

In fact, it suffices to show this theorem in the case when $X$ is a point:

\begin{theorem}[Finitary norm convergence, simpler base case]\label{main3-base-2}  
Let $F: \N \to \N$ be a function, and let $\eps > 0$.  Then there exists an integer $M^* \ge 1$ with the following property:  If $P \geq 1$, and $g: \Z_P \to [0,1]$, then there exists $1 \leq M \leq M^*$ such that
\begin{equation}\label{l2-stability-base-2}
 \| S_N g - S_{N'} g\|_{L^2(\Z_P)} \leq \eps
\end{equation}
for all $M \leq N, N' \leq F(M)$, where $S_N$ is the averaging operator $S_N g(v) := \E_{n \in [N]} g(v+n)$, and similarly for $S_{N'}$.
\end{theorem}

Indeed, Theorem \ref{main3-base} can be immediately deduced from Theorem \ref{main3-base-2} by applying the finitary Lebesgue dominated convergence theorem, Theorem \ref{fin-ldct}, using the functions
$$ f_{N,N'}(x) := \| S_N g(\cdot,x) - S_{N'} g(\cdot,x) \|_{L^2(\Z_P)}^2 \in [0,1].$$

\begin{remark} Theorem \ref{main3-base-2} is nothing more than the $l=1$ case of Theorem \ref{main2}.
\end{remark}

It remains to prove Theorem \ref{main3-base-2}.  We will not give the shortest proof of this theorem here\footnote{Indeed, one can use the Furstenberg correspondence principle to deduce Theorem \ref{main3-base-2} from the mean ergodic theorem. See also \cite{agt} for a direct proof of this theorem.}, but will instead give a more pedestrian argument which will motivate the proof of the inductive case $d > 1$ of Theorem \ref{main3} in the next section.

A crucial notion to our argument is that of an \emph{basic anti-uniform function}\footnote{Our terminology is inspired by that in \cite{gt-primes}.}.

\begin{definition}[Basic $\{1\}$-anti-uniform function]\label{baf} Let $M \geq 1$.  A \emph{basic $\{1\}$-anti-uniform function} of scale $M$ is any function $\varphi: \Z_P \to \R$ of the form
$$ \varphi(v) := \E_{n \in [M]} b(v-n)$$
for some function $b: \Z_P \to [-1,1]$.
\end{definition}

Note that any basic $\{1\}$-anti-uniform function will itself take values between $-1$ and $1$.  Furthermore, one easily verifies the Lipschitz property
\begin{equation}\label{lip}
 |\varphi(v+n) - \varphi(v)| \leq 2 \frac{|n|}{M}
\end{equation}
for all $n \in \Z$ and $v \in \Z_P$, and all basic $\{1\}$-anti-uniform functions $\varphi$ of scale $M$.  Heuristically, basic $\{1\}$-anti-uniform functions should be viewed as essentially being constant at scales below $M$.  The relevance of basic $\{1\}$-anti-uniform functions to Theorem \ref{main3-base-2} relies on the following simple lemma.  

\begin{lemma}[Lack of uniformity implies correlation with basic anti-uniform function]\label{locor}  Let $g: \Z_P \to [-1,1]$, $M \geq 1$, and $\eps > 0$ be such that
\begin{equation}\label{snog}
\|S_N g \|_{L^2(\Z_P)} \geq \eps
\end{equation}
for some $N \geq \frac{10 M}{\eps^2}$.
Then there exists a basic $\{1\}$-anti-uniform function $\varphi$ of scale $M$ such that $|\langle g, \varphi \rangle_{L^2(\Z_P)}| \geq \eps^2/2$.
\end{lemma}

\begin{proof}  We expand \eqref{snog} as
$$ \E_{v \in \Z_P} (\E_{n \in [N]} g(v+n)) (\E_{n' \in [N]} g(v+n')) \geq \eps^2.$$
Observe from the hypothesis $N \geq \frac{10M}{\eps^2}$ that
$$ |\E_{n' \in [N]} g(v+n') - \E_{n' \in [N]} \E_{m \in [M]} g(v+n'+m)| \leq \frac{\eps^2}{5}$$
for all $v \in \Z_P$,
and thus by the triangle inequality
$$ |\E_{v \in \Z_P} (\E_{n \in [N]} g(v+n)) \E_{n' \in [N]} \E_{m \in [M]} g(v+n'+m)| \geq \eps^2/2.$$
By the pigeonhole principle, we can thus find $n, n' \in [N]$ such that
$$ |\E_{v \in \Z_P} g(v+n) \E_{m \in [M]} g(v+n'+m)| \geq \eps^2/2.$$
We can rewrite this as $|\langle g, \varphi \rangle_{L^2(\Z_P)}| \geq \eps^2/2$, where
$$ \varphi(v) := \E_{m \in [M]} b(v+m)$$
and $b(v) := g(v+n'-n)$, and the claim follows.
\end{proof}

To exploit this lemma, we need to use the basic $\{1\}$-anti-uniform functions to build various factors (the finitary analogue of characteristic factors), using the construction in Lemma \ref{factor-lemma}.

We turn to the details.  Let $K \geq 1$ be the first integer larger than $\frac{10^6}{\eps^4}+2$, and $\tilde F: \N \to \N$ be a function to be chosen later (it shall depend on $F$ and $\eps$), such that $\tilde F(M) \geq M$ for all $M$.  Define the sequence 
$$ 1 \leq M_1 \leq M_2 \leq \ldots \leq M_K$$
recursively by $M_1 := 1$ and $M_{i+1} := \tilde F(M_i)$.  

By greedily iterating Lemma \ref{locor} at a rapidly diminising sequence of scales we shall obtain a useful decomposition $g = g_{U^\perp} + g_U$ where $g_{U^\perp}$ is ``low complexity'' and $g_U$ is ``negligible'' at scales between $M_{k-1}$ and $M_k$ for some $k$, in the following precise sense.

\begin{proposition}[Koopman-von Neumann type theorem]\label{kvn}  Let $g: \Z_P \to [0,1]$.  Then we can decompose $g = g_{U^\perp} + g_U$, where the two components $g_{U^\perp}, g_U: \Z_P \to [-1,1]$ have the following properties.
\begin{itemize}
\item[(i)] ($g_{U^\perp}$ anti-uniform) There exists an integer $2 \leq k \leq K$ and a basic $\{1\}$-anti-uniform function $\varphi_j$ of scale $M_j$ for eack $k \leq j \leq K$ such that $g_{U^\perp}$ is $\Y_{\geq k}$-measurable, where $\Y_{\geq k} := \Y_{\eps^2/400}(\varphi_k) \vee \ldots \vee \Y_{\eps^2/400}(\varphi_K)$, and the factors $Y_{\eps^2/400}(\varphi_j)$ are those defined in Lemma \ref{factor-lemma}.
\item[(ii)] ($g_U$ uniform) We have
\begin{equation}\label{evp}
 \| S_N g_{U} \|_{L^2(\Z_P)} \leq \eps/10
\end{equation}
for all $N \geq \frac{1000 M_{k-1}}{\eps^2}$.
\end{itemize}
\end{proposition}

\begin{remark} See \cite[Proposition 8.1]{gt-primes}, \cite[Theorem 3.9]{tao:gauss}, \cite[Theorem 4.7]{tz}, or \cite[Theorem 6.7]{gt-r4} for similar results.
\end{remark}

\begin{proof} We perform the following algorithm:

\begin{itemize}
\item Step 0.  Initialise $k=K+1$.
\item Step 1.  Set $\Y_{\geq k} := \Y_{\eps^2/400}(\varphi_k) \vee \ldots \vee \Y_{\eps^2/400}(\varphi_K)$, and then 
set $g_{U^\perp} := \E( g | \Y_{\geq k} )$ and $g_U := g - g_{U^\perp}$.  (Thus, initially, $g_{U^\perp}$ is simply the mean value $\E_{v \in \Z_P} g(v)$ of $g$.)
\item Step 2.  If \eqref{evp} holds for all $N \geq \frac{1000 M_{k-1}}{\eps^2}$ then \textbf{STOP}.  Otherwise, we apply 
Lemma \ref{locor} to locate a basic $\{1\}$-anti-uniform function $\varphi_{k-1}$ of scale $M_{k-1}$ such that $|\langle g_U, \varphi_{k-1} \rangle_{L^2(\Z_P)}| \geq \eps^2/200$.
\item Step 3.  We decrement $k$ to $k-1$.  If $k=1$ then we \textbf{STOP} with an error; otherwise we return to Step 1.
\end{itemize}

If this algorithm terminates at some $k \geq 2$ then we are done, so suppose instead for contradiction that the algorithm goes all the way down to $k=1$.  Then we have constructed $\varphi_1,\ldots,\varphi_K$ such that
$$ \left|\left\langle g - \E(g | \Y_{\geq j+1} ), \varphi_j \right\rangle_{L^2(\Z_P)}\right| \geq \eps^2/200$$
for all $1 \leq j \leq K$.  On the other hand, by Lemma \ref{factor-lemma}(i) we have
$$ \left\| \varphi_j  - \E(\varphi_j | \Y_{\geq j}) \right\|_{L^\infty} \leq \eps^2/400$$
and hence by the triangle inequality (and the fact that $g$ takes values in $[0,1]$) we have
$$ \left|\left\langle g - \E(g | \Y_{\geq j+1} ), \E(\varphi_j | \Y_{\geq j} ) \right\rangle_{L^2(\Z_P)}\right| \geq \eps^2/400.$$
We can rewrite the left-hand side as
$$ \left|\left\langle \E(g|\Y_{\geq j} )
- \E(g | \Y_{\geq j+1} ), \E(\varphi_j | \Y_{\geq j} ) \right\rangle_{L^2(\Z_P)}\right|$$
and thus by Cauchy-Schwarz
$$ \left\| \E(g|\Y_{\geq j} )
- \E(g | \Y_{\geq j+1} ) \right\|_{L^2(\Z_P)} \geq \eps^2/400$$
and thus by Pythagoras' theorem
$$ \left\| \E(g|\Y_{\geq j} ) \right\|_{L^2(\Z_P)}^2
\geq \left\| \E(g|\Y_{\geq j+1} ) \right\|_{L^2(\Z_P)}^2 + \frac{\eps^4}{10^6}$$
for all $1 \leq j \leq K$.  On the other hand, the quantities $\| \E(g|\Y_{\geq j} ) \|_{L^2(\Z_P)}^2$ clearly range between $0$ and $1$.  These facts contradict the definition of $K$.  The claim follows.
\end{proof}

We apply this proposition to obtain $2 \leq k \leq K$, basic $\{1\}$-anti-uniform functions $\varphi_k,\ldots,\varphi_K$, and a decomposition $g = g_{U^\perp} + g_U$ with the stated properties.

Let $M$ be the first integer greater than $\frac{1000 M_{k-1}}{\eps}$, and let $M \leq N, N' \leq F(M)$.  To prove Theorem \ref{main3-base-2}, it will suffice to show that, for $\tilde F$ chosen appropriately depending on $F$ and $\eps$, 
$$ \| S_N g - S_{N'} g \|_{L^2(\Z_P)} \leq \eps,$$
since $M$ will be bounded by some quantity $M_*$ depending on $\eps$ and $\tilde F$, and thus ultimately on $F$ and $\eps$.  From \eqref{evp} we already have
$$ \| S_N g_U \|_{L^2(\Z_P)}, \| S_{N'} g_U \|_{L^2(\Z_P)} \leq \eps/10$$
so by the triangle inequality it will suffice to show that
\begin{equation}\label{guv}
\| S_N g_{U^\perp} - S_{N'} g_{U^\perp} \|_{L^2(\Z_P)} \leq \eps/10.
\end{equation}

Now the function $g_{U^\perp}$ takes values between $0$ and $1$, and is measurable with respect to the factor $\Y_{\geq k}$.  From Lemma \ref{factor-lemma}, this factor has $O_{K,\eps}(1) = O_\eps(1)$ atoms, each of which is the intersection of atoms coming from the individual factors $\Y_{\eps^2/400}(\varphi_k), \ldots, \Y_{\eps^2/400}(\varphi_K)$.  Applying Lemma \ref{factor-lemma} repeatedly, we thus see for every $\eta_1 > 0$ there exists a polynomial $\Psi: \R^{K-k+1} \to \R$ of $K-k+1$ variables with degree and coefficients $O_{K,\eps,\eta_1}(1) = O_{\eps,\eta_1}(1)$ such that
$$ \| g_{U^\perp} - \Psi( \varphi_k, \ldots, \varphi_K ) \|_{L^1(\Z_P)} \ll_{\eps} \eta_1$$
and
$$ \| g_{U^\perp} - \Psi( \varphi_k, \ldots, \varphi_K ) \|_{L^\infty(\Z_P)} \ll_{\eps} 1.$$
By H\"older's inequality we conclude that
$$ \| g_{U^\perp} - \Psi( \varphi_k, \ldots, \varphi_K ) \|_{L^2(\Z_P)}^2 \ll_{\eps} \eta_1;$$
since $S_N$ is a contraction on $L^2$, we conclude that
$$ \| S_N g_{U^\perp} - S_N \Psi( \varphi_k, \ldots, \varphi_K ) \|_{L^2(\Z_P)}^2 \ll_{\eps} \eta_1;$$
Thus, if we choose $\eta_1$ sufficiently small depending on $\eps$, we see from the triangle inequality that \eqref{guv} will follow if we can show
\begin{equation}\label{guv-2} 
 \| S_N \Psi( \varphi_k, \ldots, \varphi_K ) - S_{N'} \Psi( \varphi_k, \ldots, \varphi_K ) \|_{L^2(\Z_P)} \leq \eps/20.
\end{equation}
We now fix $\eta_1 = \eta_1(\eps)$ so that the above argument is valid. From \eqref{lip} (and the monotonicity of the $M_j$) we have
$$ \varphi_j(v+n) = \varphi_j(v) + O\left( \frac{F(M)}{M_k} \right)$$
for all $k \leq j \leq K$ and $n \in [N] \cup [N']$.  By the bounds on $\Psi$ (and the fact that the $\varphi_j$ have magnitude $O(1)$) we conclude that
\begin{align*}
\Psi( \varphi_k, \ldots, \varphi_K )(v+n) &= 
\Psi( \varphi_k, \ldots, \varphi_K )(v) + O_{K,\eps,\eta_1}\left( \frac{F(M)}{M_k} \right) \\
&= \Psi( \varphi_k, \ldots, \varphi_K )(v) + O_{\eps}\left( \frac{F(M)}{\tilde F( M_{k-1} )} \right);
\end{align*}
averaging in $n$, we obtain
$$ S_N \Psi( \varphi_k, \ldots, \varphi_K ), S_{N'} \Psi( \varphi_k, \ldots, \varphi_K )
= \Psi( \varphi_k, \ldots, \varphi_K ) + O_{\eps}\left( \frac{F(M)}{\tilde F( M_{k-1} )} \right).$$
Thus we can bound the left-hand side of \eqref{guv-2} by $O_{\eps}\left( \frac{F(M)}{\tilde F( M_{k-1} )} \right)$.  If we then choose $\tilde F$ to grow sufficiently quickly depending on $F$ and $\eps$ we obtain the desired claim (setting $M^* := M_K$).
This concludes the proof of Theorem \ref{main3-base-2}, and hence the $d=1$ case of Theorem \ref{main3}.

\section{The inductive case}\label{inductive-sec}

To complete the proof of Theorem \ref{main3} (and thus Theorem \ref{main}) it remains to verify the inductive step of Theorem \ref{main3}, i.e. to deduce Theorem \ref{main3} for some fixed $d > 1$ assuming inductively that this theorem has already been established for all smaller values of $d$.  Fortunately it turns out that the arguments of the preceding section extend without much difficulty to handle this case.  The one twist will be that the basic anti-uniform functions will have higher complexity (they are averages of complexity $d-1$), and in particular will not obey the simple Lipschitz property \eqref{lip}; however, they will be approximable by functions of complexity $d-1$ or less and will thus be treatable by the induction hypothesis\footnote{In ergodic theory terminology, the complexity $d$ case (with $J=1$) is being viewed as a kind of ``weakly mixing extension'' of the complexity $d-1$ case (with $J > 1$), with the latter serving as a kind of ``characteristic factor'' for the former.  Similarly, the $J > 1$ case at a given complexity is a kind of ``finite rank extension'' of the $J=1$ case.}.

Before we begin the rigorous argument, let us give an informal discussion to try to motivate the strategy of proof.  For simplicity let us just discuss the case $d=2$ and $l=3$, with $X$ equal to a point, and consider the convergence of averages $\Delta_N(f)$, where $f$ has complexity at most $(2,1)$, and specifically $f$ takes the form
$$ f(v_1,v_2,v_3) = g_{\{1,2\}}(v_1,v_2) g_{\{2,3\}}(v_2,v_3) g_{\{3,1\}}(v_3,v_1)$$
for some functions\footnote{To be completely consistent with our other notation, we should actually make each of $g_{\{1,2\}}, g_{\{2,3\}}, g_{\{3,1\}}$ equal to a function on $\Z_P^3$ which is constant in one of the variables $v_1,v_2,v_3$, but we will not do so here to simplify the formulas slightly.} $g_{\{1,2\}}, g_{\{2,3\}}, g_{\{3,1\}}: \Z_P^2 \to [-1,1]$.
Then the average $\Delta_N(f)$ can be written explicitly as
$$ \Delta_N(f)(v_1,v_2) = \E_{n \in [N]} g_{\{1,2\}}(v_1,v_2) g_{\{2,3\}}(v_2,-v_1-v_2-n) g_{\{3,1\}}(-v_1-v_2-n,v_1).$$
The $g_{\{1,2\}}(v_1,v_2)$ factor comes out of the average (cf. \eqref{gih}) and is therefore uninteresting.  We shall thus assume $g_{\{1,2\}} \equiv 1$ and so
\begin{equation}\label{deltan}
 \Delta_N(f)(v_1,v_2) = \E_{n \in [N]} g_{\{2,3\}}(v_2,-v_1-v_2-n) g_{\{3,1\}}(-v_1-v_2-n,v_1).
 \end{equation}
Now suppose that we are in the ``compact'' or ``finite rank'' case in which $g_{\{2,3\}}$ and $g_{\{3,1\}}$ were actually complexity $1$ objects, for instance suppose we had
$$ g_{\{2,3\}}(v_2,v_3) = h_2(v_2) h_3(v_3) \hbox{ and } g_{\{3,1\}}(v_3,v_1) = k_3(v_3) k_1(v_1)$$
for some functions $h_2,h_3,k_3,k_1: \Z_P \to [-1,1]$.  Then the average simplifies to
$$ \Delta_N(f)(v_1,v_2) = h_2(v_2) k_1(v_1) \E_{n \in [N]} h_3 k_3(-v_1-v_2-n).$$
The convergence of this average can then be easily deduced from the $d=1$ theory of the previous section.  Similarly we expect to be able to handle the case when $g_{\{2,3\}}$ and $g_{\{3,1\}}$ are of complexity $(1,J)$ for some bounded $J$, i.e. they are a bounded combination of tensor products of functions of one variable.

Now let us consider the opposing case in which $g_{\{2,3\}}$ (say) does \emph{not} behave at all like a tensor product of one variable, so much so that they behave ``orthogonally'' to any such tensor products.  A little more precisely, let us suppose that correlations of the form
\begin{equation}\label{evw}
\E_{v_2 \in w_2 + [N']; v_3 \in w_3 + [N']} g_{\{2,3\}}(v_2,v_3) h_2(v_2) h_3(v_3)
\end{equation}
are always small for ``generic'' base points $w_2, w_3 \in \Z_P$ and arbitrary bounded functions $h_2, h_3: \Z_P \to [-1,1]$ (we will not attempt to make these assertions rigorous here), and for various values of $N'$ which we shall leave vague here.  In that ``weakly mixing'' case, it turns out that the averages $\Delta_N(f)$ are in fact quite small in norm.  To see this, let us write 
$$ \|\Delta_N(f)\|_{L^2}^2 = \frac{1}{P^2} \sum_{v_1,v_2 \in \Z_P} \Delta_N(f)(v_1,v_2)
\E_{n \in [N]} g_{\{2,3\}}(v_2,-v_1-v_2-n) g_{\{3,1\}}(-v_1-v_2-n,v_1)$$
and then rewrite the right-hand side as
$$ \frac{1}{NP^2} \sum_{v_1 \in \Z_P} \sum_{v_2,v_3 \in\Z_P: -\Sigma(v_1,v_2,v_3) \in [N]} 
g_{\{2,3\}}(v_2,v_3) \Delta_N(f)(v_1,v_2) g_{\{3,1\}}(v_3,v_1).$$
But observe that for any fixed $v_1$, the inner sum is (up to some normalising factors) the correlation between $g_{\{2,3\}}(v_2,v_3)$ and a tensor product of functions of $v_2$ and $v_3$ separately.  This sum is over a diagonal region $\{ (v_2,v_3): -\Sigma(v_1,v_2,v_3) \in [N]\}$, but we can approximately split this region into squares of length $N'$ for some $N'$ a bit smaller than $N$ and use the smallness of \eqref{evw} to then conclude that $\Delta_N(f)$ is small in $L^2$.

To summarise so far, we have given heuristics to justify some sort of convergence in the extreme cases when both $g_{\{2,3\}}$ and $g_{\{3,1\}}$ are ``compact'', and when at least one of $g_{\{2,3\}}$ and $g_{\{3,1\}}$ are ``weakly mixing''.  The rest of the proof then hinges on a Koopman-von Neumann type structure theorem (as in the previous section) that allows us to split arbitrary functions into compact and weakly mixing components, allowing us to deduce the general case from these two special cases.

We turn to the details.  Fix $d > 1$, and assume inductively that Theorem \ref{main3} has already been established for all smaller values of $d$.  We allow all implied constants to depend on $l$ and $d$.  By increasing $F$ if necessary we may assume that $F(M) \geq M$ for all $M$.

We can first repeat several of the reductions already employed in the previous section.  For instance, we can quickly reduce to the case $M_*=J=1$ by using exactly the same arguments used in the $d=1$ case.  Similarly, by using Theorem \ref{fin-ldct} as before we may reduce $X$ to a point.  If we write $g = \prod_{e \subset \{1,\ldots,l+1\}: |e|=d} g_e$, where $g_e: \Z_P^{l+1} \to [-1,1]$ is $e$-measurable, then as before the contribution of those $e$ for which $e \subset \{1,\ldots,l\}$ can be absorbed using the module identity \eqref{gih}.  Our task is now to establish the following.

\begin{theorem}[Finitary norm convergence, inductive step]\label{main3b}  Let $1 < d \leq l$, and suppose that Theorem \ref{main3} has already been established for smaller values of $d$.  Let $\I$ be the collection of all subsets $e$ of $\{1,\ldots,l+1\}$ such that $|e| = d$ and $l+1 \in e$.  Let $F: \N \to \N$ be a function, and let $\eps > 0$.  Then there exists an integer $M^* \ge 1$ with the following property:  If $P \geq 1$, and $g_e: \Z_P^{l+1} \to [-1,1]$ is $e$-measurable for all $e \in \I$, then there exists an integer $1 \leq M \leq M^*$ such that
\begin{equation}\label{l2-stability-induct}
\left\| \Delta_N\left(\prod_{e \in \I} g_e\right) - \Delta_{N'}\left(\prod_{e \in \I} g_e\right) \right\|_{L^2(\Z_P^{l})} \leq \eps
\end{equation}
for all $M \leq N, N' \leq F(M)$.  
\end{theorem}

As in the previous section, a key concept will be that of an anti-uniform function, although now this function will be adapted to the index set $e$.

\begin{definition}[Basic $e$-anti-uniform function]\label{bafe} Let $M \geq 1$, and let $e \in \I$.  A \emph{basic $e$-anti-uniform function} of scale $M$ is any function $\varphi_e: \Z_P^{l+1} \to \R$ of the form
$$ \varphi_e(v) := \E_{m \in [M]} \prod_{i \in e} b_i\left( v_{e \backslash \{i\}}, \Sigma(v_e) + m \right)$$
where for each $i \in e$, $b_i: \Z_P^{e \backslash \{i\}} \times \Z_P \to [-1,1]$ is a function, and for each $v = (v_1,\ldots,v_{l+1}) \in \Z_P^{l+1}$, $v_e := (v_j)_{j \in e} \in \Z_P^e$ and $v_{e \backslash \{i\}} := (v_j)_{j \in e \backslash \{i\}} \in \Z_P^{e \backslash \{i\}}$ are projections of $v$.
\end{definition}

Observe that this definition generalises Definition \ref{baf}, which considered the case $l=0$ and $e=\{1\}$.  Also note that any basic $e$-anti-uniform function $\varphi_e$ of scale $M$ is going to be $e$-measurable and take values in $[-1,1]$. 

\begin{example} If $l=2$ and $e = \{1,2\}$, and $b_1, b_2: \Z_P^2 \to [-1,1]$, then any function of the form
$$ \varphi_e(v_1,v_2,v_3) = \E_{m \in [M]} b_1( v_2, v_1+v_2+m ) b_2( v_1, v_1+v_2+m )$$
is a basic $e$-anti-uniform of scale $M$.  
\end{example}

We have a generalisation of Lemma \ref{locor}:

\begin{lemma}[Lack of uniformity implies correlation with basic anti-uniform function]\label{locor-gen}  Let $M \geq 1$ and $\eps > 0$.  For each $e \in \I$, let $g_e: \Z_P^{l+1} \to [-1,1]$ be an $e$-measurable function, and suppose that
\begin{equation}\label{veg}
 \left\| \Delta_N\left(\prod_{e \in \I} g_e\right) \right\|_{L^2(\Z_P^{l})} \geq \eps
\end{equation}
for some $N \geq \frac{10M}{\eps^2}$.  Then for every $e_0 \in \I$, there exists a basic $e_0$-anti-uniform function $\varphi_{e_0}$ such that $|\langle g_{e_0}, \varphi_{e_0}\rangle_{L^2(\Z_P^{l+1})}| \geq \eps^2/2$.
\end{lemma}

\begin{proof}
From Definition \ref{dap}, we have
$$ \Delta_N\left(\prod_{e \in \I} g_e\right)( v_1,\ldots,v_l ) := \frac{1}{N} \sum_{v_{l+1}: -\Sigma(v) \in [N]} \prod_{e \in \I} g_e(v)$$
where $v = (v_1,\ldots,v_{l+1})$.  Squaring \eqref{veg}, we obtain
$$ \sum_{(v_1,\ldots,v_l) \in \Z_P^l} \Delta_N\left(\prod_{e \in \I} g_e\right)(v_1,\ldots,v_l)
\sum_{v_{l+1}: -\Sigma(v) \in [N]} \prod_{e \in \I} g_e(v) \geq \eps^2 N P^l;$$
if we then let $h: \Z_P^{l+1} \to [-1,1]$ be the function
$$ h(v_1,\ldots,v_{l+1}) := \Delta_N\left(\prod_{e \in \I} g_e\right)(v_1,\ldots,v_l),$$
we then obtain
$$ \sum_{v \in \Z_P^{l+1}: -\Sigma(v) \in [N]} h(v) \prod_{e \in \I} g_e(v) \geq \eps^2 N P^l.$$
Observe that for each $e \in \I \backslash \{e_0\}$, $g_e$ will be $\{1,\ldots,l+1\} \backslash \{i\}$-measurable for some $i \in e_0$.  The function $h$ obeys the same property; indeed, $h$ is clearly $\{1,\ldots,l+1\} \backslash \{l+1\}$-measurable, and $l+1$ lies in $e_0$ by definition of $\I$.
$$ h(v) \prod_{e \in \I} g_e(v) = g_{e_0}(v) \prod_{i \in e_0} b_i(v)$$
where $b_i: \Z_P^{l+1} \to [-1,1]$ is a $\{1,\ldots,l+1\} \backslash \{i\}$-measurable function.  Thus we have
$$ \sum_{v_{e_0} \in \Z_P^{e_0}} g_{e_0}(v_{e_0})
\sum_{v_{e^c_0} \in \Z_P^{e^c_0}: -\Sigma(v_{e_0}) - \Sigma(v_{e^c_0}) \in [N]} \prod_{i \in e_0} b_i(v_{e_0}, v_{e^c_0}) \geq \eps^2 N P^l,$$
where $e^c_0 := \{1,\ldots,l+1\} \backslash e_0$, and we abuse notation by identifying the $e_0$-measurable function $g_{e_0}$ with a function on $\Z_P^{e_0}$.  Since $e^c_0$ has cardinality $l+1-d > 0$, we can write $e^c_0 = \{j\} \cup f$ for some $j \in \{1,\ldots,l+1\}$ and some $f \subset \{1,\ldots,l+1\}$ of cardinality $l-d$.  By the pigeonhole principle, we may thus find $v_f \in \Z_P^f$ such that
$$ \sum_{v_{e_0} \in \Z_P^{e_0}} g_{e_0}(v_{e_0})
\sum_{v_j \in \Z_P^{e^c_0}: -\Sigma(v_{e_0}) - \Sigma(v_f) - v_j \in [N]} \prod_{i \in e_0} b_i(v_{e_0}, v_j, v_f) \geq \eps^2 N P^d.$$
Fix this $v_f$.  Since $N \geq \frac{10M}{\eps^2}$, we can shift $[N]$ by $m$ for any $m \in [M]$ and only pick up an error of at most $\eps^2 N P^d/2$, thus
$$ \sum_{v_{e_0} \in \Z_P^{e_0}} g_{e_0}(v)
\sum_{v_j \in \Z_P^{e^c_0}: -\Sigma(v_{e_0}) - \Sigma(v_f) - v_j \in [N]+m} \prod_{i \in e_0} b_i(v_{e_0}, v_j, v_f) \geq \eps^2 N P^d/2$$
for all $m \in [M]$.  Summing this over all $m \in [M]$ we obtain
$$ \sum_{n \in [N]} \sum_{v_{e_0} \in \Z_P^{e_0}} g_{e_0}(v)
\sum_{v_j \in \Z_P^{e^c_0}: -\Sigma(v_{e_0}) - \Sigma(v_f) - v_j \in n+[M]} \prod_{i \in e_0} b_i(v_{e_0}, v_j, v_f) \geq \eps^2 NM P^d/2.$$
By the pigeonhole principle we may thus find $n \in [N]$ such that
$$ \sum_{v_{e_0} \in \Z_P^{e_0}} g_{e_0}(v)
\sum_{v_j \in \Z_P^{e^c_0}: -\Sigma(v_{e_0}) - \Sigma(v_f) - v_j \in n+[M]} \prod_{i \in e_0} b_i(v_{e_0}, v_j, v_f) \geq \eps^2 M P^d/2.$$
If we define $\tilde b_i: \Z_P^{e_0} \times \Z_P \to [-1,1]$ to be the function
$$ \tilde b_i(v_{e_0}, w) := b_i( v_{e_0}, -\Sigma(v_f) - w - n, v_f )$$
then we have
$$ \sum_{v_{e_0} \in \Z_P^{e_0}} g_{e_0}(v_{e_0}) \sum_{m \in [M]} \prod_{i \in e_0} \tilde b_i( \tilde v_{e_0}, \Sigma(v_{e_0}) + m ) \geq \eps^2 M P^d/2,$$
or in other words
$$ \E_{v_{e_0} \in \Z_P^{e_0}} g_{e_0}(v_{e_0}) \E_{m \in [M]} \prod_{i \in e_0} \tilde b_i( \tilde v_{e_0}, \Sigma(v_{e_0}) + m ) \geq \eps^2/2.$$
If we now add some dummy variables $v_k$ for all $k \in e_0^c$, we obtain the claim.
\end{proof}

Now let $K \geq 1$ be the first integer larger than $\frac{10^6|\I|^5}{\eps^4}+2$, and $\tilde F: \N \to \N$ be a function to be chosen later (it shall depend on $F$ and $\eps$), such that $\tilde F(M) \geq M$ for all $M$.  Once again, we define the sequence 
$$ 1 \leq M_1 \leq M_2 \leq \ldots \leq M_K$$
recursively by $M_1 := 1$ and $M_{i+1} := \tilde F(M_i)$.  By running the proof of Proposition \ref{kvn} ``in parallel'' for each of the $g_e$ simultaneously, we now show

\begin{proposition}[Koopman-von Neumann type theorem]\label{kvn-2}  For each $e \in \I$, let $g_e: \Z_P^{l+1} \to [0,1]$ be an $e$-measurable function.  Then there exists $2 \leq k \leq K+1$ and decompositions $g_e = g_{e,U^\perp} + g_{e,U}$ for all $e \in \I$, where $g_{e,U^\perp}, g_{e,U}: \Z_P^{l+1} \to [-1,1]$ are $e$-measurable functions with the following properties.
\begin{itemize}
\item[(i)] ($g_{e,U^\perp}$ anti-uniform) For each $e \in \I$, there exists a basic $e$-anti-uniform function $\varphi_{e,j}$ of scale $M_j$ for each $k \leq j \leq K$ such that $g_{e,U^\perp}$ is $\Y_{e,\geq k}$-measurable, where $\Y_{e,\geq k} := \Y_{\eps^2/(400|\I|^2)}(\varphi_{e,k}) \vee \ldots \vee \Y_{\eps^2/(400|\I|^2)}(\varphi_{e,K})$.
\item[(ii)] ($g_{e,U}$ uniform) For any $e \in \I$, we have
\begin{equation}\label{evp2}
 \| \Delta_N(g_{e,U} \prod_{e' \in \I \backslash \{e\}} h_{e'}) \|_{L^2(\Z_P^{l})} \leq \frac{\eps}{10|\I|}
\end{equation}
for all $N \geq \frac{1000 |\I|^2 M_{k-1}}{\eps^2}$ and all $e'$-measurable $h_{e'}: \Z_P^{l+1} \to [-1,1]$ for $e' \in \I \backslash \{e\}$.
\end{itemize}
\end{proposition}

\begin{remark} This result is a ``weak hypergraph regularity lemma'', akin to the ``weak regularity lemma'' of Frieze and Kannan \cite{fkannan}.  One can also develop stronger regularity lemmas (in which one obtains \emph{local} regularity and not just \emph{global} regularity), similar for instance to those in \cite{tao:hyper}, by replacing the ``single-loop'' greedy algorithm argument presented here by a ``double-loop'' one, but they will not be necessary for our purposes here.
\end{remark}

\begin{proof} The argument shall closely follow the proof of Proposition \ref{kvn}.
We perform the following algorithm:

\begin{itemize}
\item Step 0.  Initialise $k=K+1$.
\item Step 1.  For each $e \in \I$, set $\Y_{e,\geq k} := \Y_{\eps^2/(400|\I|^2)}(\varphi_{e,k}) \vee \ldots \vee \Y_{\eps^2/(400|\I|^2)}(\varphi_{e,K})$,
$g_{e,U^\perp} := \E( g_e | \Y_{e,\geq k} )$ and $g_{e,U} := g_e - g_{e,U^\perp}$.  
\item Step 2.  If \eqref{evp2} holds for all $N \geq \frac{1000|\I|^2 M_{k-1}}{\eps^2}$, all $e \in \I$, and all $e'$-measurable $h_{e'}: \Z_P^{l+1} \to [-1,1]$ then \textbf{STOP}.  Otherwise, we apply 
Lemma \ref{locor-gen} to locate an $e \in \I$ and a basic $e$-anti-uniform function $\varphi_{e,k-1}$ and scale $M_{k-1}$ such that $|\langle g_{e,U}, \varphi_{e,k-1} \rangle_{L^2(\Z_P^{l+1})}| \geq \eps^2/(200|\I|^2)$.  For all the $e'$ in $\I$ that are not equal to $e$, we set $\varphi_{e',k-1}$ to be an arbitrary basic $e'$-anti-uniform function of scale $M_{k-1}$ (e.g. one could set $\varphi_{e',k-1} := 1$).
\item Step 3.  We decrement $k$ to $k-1$.  If $k=1$ then we \textbf{STOP} with an error; otherwise we return to Step 1.
\end{itemize}

Once again, we are done if this algorithm terminates at some $k \geq 2$, so suppose instead for contradiction that the algorithm goes all the way down to $k=1$.  Then, by construction, we have constructed $\varphi_{e,j}$ for $e \in \I$ and $1 \leq j \leq K$, with the property that for every $1 \leq j \leq K$ there exists $e \in \I$ such that
$$ \left|\left\langle g_e - \E(g_e | \Y_{e, \geq j+1} ), \varphi_{e,j} \right\rangle_{L^2(\Z_P^{l+1})}\right| \geq \frac{\eps^2}{200|\I|^2}.$$
By arguing exactly as in the proof of Proposition \ref{kvn}, we then conclude that
$$ \left\| \E(g_e|\Y_{e, \geq j} ) \right\|_{L^2(\Z_P^{l+1})}^2
\geq \left\| \E(g_e|\Y_{e, \geq j+1} ) \right\|_{L^2(\Z_P^{l+1})}^2 + \frac{\eps^4}{10^6|\I|^4}$$
for this value of $e$.  On the other hand, from Pythagoras' theorem we have
$$ \left\| \E(g_{e'}|\Y_{e', \geq j} ) \right\|_{L^2(\Z_P^{l+1})}^2
\geq \left\| \E(g_{e'}|\Y_{e', \geq j+1} ) \right\|_{L^2(\Z_P^{l+1})}^2$$
for all other values of $e' \in \I$.  Thus if we define
$$ c_j := \sum_{e' \in \I} \left\| \E(g_{e'}|\Y_{e', \geq j} ) \right\|_{L^2(\Z_P^{l+1})}^2$$
then we have
$$ c_j \geq c_{j+1} + \frac{\eps^4}{10^6 |\I|^4}.$$
On the other hand, $c_j$ varies between $0$ and $|\I|$.  This contradicts the choice of $K$, and Proposition \ref{kvn-2} follows.
\end{proof}

We apply this proposition to obtain $2 \leq k \leq K$, basic $e$-anti-uniform functions $\varphi_{e,j}$ for $e \in \I$ and $k \leq j \leq K$, and decompositions $g_e = g_{e,U^\perp} + g_{e,U}$ with the stated properties.

Let $M_{**}$ be the first integer greater than $\frac{1000|\I|^2 M_{k-1}}{\eps^2}$, and let $M^{**}$ be the first integer such that $F(M^{**}) \geq M_k^{1/4}$ (so in particular $M^{**} \leq M_k^{1/4}+1$).  Thus
$$ 1 \leq M_{k-1} \leq M_{**} \leq M^{**} \leq M_k \leq \ldots \leq M_K.$$
To prove Theorem \ref{main3b} (with $M^* := M_K$), it will suffice to show that, for $\tilde F$ chosen appropriately depending on $F$ and $\eps$, that there exists $M_{**} \leq M < M^{**}$ such that
\begin{equation}\label{vegan}
\left\| \Delta_N\left(\prod_{e \in \I} g_e\right) - \Delta_{N'}\left(\prod_{e \in \I} g_e\right) \right\|_{L^2(\Z_P^{l})} \leq \eps
\end{equation}
for all $M \leq N,N' \leq F(M)$.  Note that since $M_k = \tilde F(M_{k-1})$, we can make $M^{**}$ larger than any specified function of $M_*$ by choosing $\tilde F$ to be sufficiently rapidly growing.

Let us enumerate $\I$ arbitrarily as $\I = \{e_1,\ldots,e_{|\I|}\}$.  From \eqref{evp2} we have
$$  
 \left\| \Delta_N\left(g_{e_j,U} \left(\prod_{1 \leq j' < j} g_{e_{j'},U^\perp}\right) \left(\prod_{j < j' \leq |\I|} g_{e_{j''}}\right)\right)\right\|_{L^2(\Z_P^{l})} \leq \eps/(10|\I|)$$
for all $1 \leq j \leq |\I|$ and all $N \geq M_*$.  From the standard telescoping identity
$$ \prod_{j=1}^{|\I|} g_{e_j} - \prod_{j=1}^{|\I|} g_{e_j,U^\perp}
 = \sum_{j=1}^{|\I|} g_{e_j,U} \left(\prod_{1 \leq j' < j} g_{e_{j'},U^\perp}\right) \left(\prod_{j < j' \leq |\I|} g_{e_{j''}}\right)$$
and the triangle inequality, we conclude that
$$ \left\| \Delta_N\left(\prod_{e \in \I} g_e\right) - \Delta_{N}\left(\prod_{e \in \I} g_{e,U^\perp}\right)\right\|_{L^2(\Z_P^{l})} \leq \eps/10.$$
By the triangle inequality again, we see that to show \eqref{vegan}, it suffices to find $M_{**} \leq M < M^{**}$ such that
\begin{equation}\label{delgan}
\left\| \Delta_N\left(\prod_{e \in \I} g_{e,U^\perp}\right) - \Delta_{N'}\left(\prod_{e \in \I} g_{e,U^\perp}\right)\right\|_{L^2(\Z_P^{l})}  \leq \eps/10
\end{equation}
for all $M \leq N,N' \leq F(M)$.

We have reduced to the ``characteristic factor'' of anti-uniform functions, and will now break these functions up into their basic components.  Let $\eta_1 > 0$ be a small quantity to be chosen later.  By applying Lemma \ref{factor-lemma} precisely as in the preceding section, we see that for every $e \in \I$ there exists a polynomial $\Psi_e: \R^{K-k+1} \to \R$ of $K-k+1$ variables with degree and coefficients $O_{\eps,\eta_1}(1)$ such that
$$ \| g_{e,U^\perp} - \Psi_e( \varphi_{e,k}, \ldots, \varphi_{e,K} ) \|_{L^1(\Z_P^{l+1})} \ll_{\eps} \eta_1$$
and
$$ \| g_{e,U^\perp} - \Psi_e( \varphi_{e,k}, \ldots, \varphi_{e,K} ) \|_{L^\infty(\Z_P^{l+1})} \ll_{\eps} 1$$
(note that as we are allowing implied constants to depend on $l$ and $d$, we have $|\I| = O(1)$ and $K = O_\eps(1)$).
Now, for any $e$-measurable function $h_e: \Z_P^{l+1} \to \R$ and any integer $n$, one can use the $e$-measurability to check that the function $(v,v_{l+1}) \mapsto h_e( (v, -\Sigma(v)-n) )$ is a permutation of $h_e$ and thus has the same $L^1$ norm.  Averaging this in $n$ using Minkowski's inequality we conclude that
$$ \| \Delta_N( h_e ) \|_{L^1(\Z_P^{l})} \leq \| h_e \|_{L^1(\Z_P^{l+1})}$$
for any $e$-measurable function $h_e: \Z_P^{l+1} \to \R$, and thus
$$ \| \Delta_N( h_e b) \|_{L^1(\Z_P^{l})} \leq \| h_e \|_{L^1(\Z_P^{l+1})}$$
for any $e$-measurable function $h_e: \Z_P^{l+1} \to \R$ and any function $b: \Z_P^{l+1} \to [-1,1]$.  Because of this and many applications of the triangle inequality we see that
$$
\left\| \Delta_N\left(\prod_{e \in \I} g_{e,U^\perp}\right) - \Delta_N(h) \right\|_{L^1(\Z_P^{l})} \ll_\eps \eta_1$$
and
$$
\left\| \Delta_N\left(\prod_{e \in \I} g_{e,U^\perp}\right) - \Delta_N(h) \right\|_{L^\infty(\Z_P^{l})} \ll_\eps 1$$
where
\begin{equation}\label{hek}
 h := \prod_{e \in \I} \Psi_e( \varphi_{e,k}, \ldots, \varphi_{e,K} ).
\end{equation}
In particular, we have
$$
\left\| \Delta_N\left(\prod_{e \in \I} g_{e,U^\perp}\right) - \Delta_N(h) \right\|_{L^2(\Z_P^{l})} \ll_\eps \eta_1^{1/2}.$$
Similarly for $N$ replaced by $N'$. Thus if we choose $\eta_1$ sufficiently small depending on $\eps$, we see from the triangle inequality that to show \eqref{delgan} it suffices to show that there exists $M_{**} \leq M \leq M^{**}$ such that
\begin{equation}\label{hon}
\| \Delta_N(h) - \Delta_{N'}(h)\|_{L^2(\Z_P^l)} \leq \eps/20,
\end{equation}
for all $M \leq N, N' \leq F(M)$.

Henceforth we fix $\eta_1$ depending on $\eps$ so that the above reductions hold.

In principle, the induction hypothesis should now let us conclude the argument.  Unfortunately, the function $h$ is not quite a function of complexity $d-1$, because of the localisations to scale $M_j$ present inside the basic $e$-anti-uniform functions $\varphi_{e,j}$.  Fortunately (as in the previous section), these scales are very large, indeed $M_j \geq M_k \geq \tilde F(M_{k-1})$, and since we have the freedom to choose $\tilde F$ at will, this localisation will end up causing no difficulty.

We turn to the details.  It will be convenient to localise the spatial variable to the scale $L := \lfloor M_k^{1/2} \rfloor$; note that this scale is intermediate between the coarse scales $M_k,\ldots,M_K$ and the fine scales $M_{**}, M^{**}$.  We can rewrite the left-hand side of \eqref{hon} as
$$ \left( \E_{v \in \Z_P^l, w \in [L]^l} |\Delta_N(h)(v+w) - \Delta_{N'}(h)(v+w)|^2 \right)^{1/2}$$
which we expand a little further using Definition \ref{dap} as
\begin{equation}\label{evil}
\left( \E_{v \in \Z_P^l, w \in [L]^l} \left|\E_{n \in [N]} h(v+w,-\Sigma(v+w)-n) - \E_{n \in [N']} h(v+w,-\Sigma(v+w)-n)\right|^2 \right)^{1/2}.
 \end{equation}
We can approximate $h$ as an average of complexity $d-1$ functions:

\begin{lemma}[$h$ essentially has complexity $d-1$]  For $v \in \Z_P^l$, $w \in [L]^l$, and $n \in [N] \cup [N']$, we can write 
$$ h(v+w,-\Sigma(v+w)-n) = \E_{\vec m \in \M} \varphi_{v,\vec m}(w, -\Sigma(w)-n) + O_\eps( M_k^{-1/2} ),$$
where $\M$ is a finite set, and for each $\vec m \in \M$, $\varphi_{v,\vec m}: \Z^{l+1} \times \Z \to \R$ is an elementary function of complexity at most $(d-1,O_\eps(1))$.
\end{lemma}

\begin{remark} The parameter $\vec m \in \M$ shall play a ``passive'' role and will eventually be absorbed into a probability space $X$ when we apply the induction hypothesis.
\end{remark}

\begin{proof} From \eqref{hek} we know that $h(v+w,-\Sigma(v+w)-n)$ is a polynomial combination of the quantities $\varphi_{e,j}(v+w,-\Sigma(v+w)-n)$.  On the other hand, from Definition \ref{bafe} we can write
$$ \varphi_{e,j}(v+w,-\Sigma(v+w)-n) := \E_{m_j \in [M_j]} \prod_{i \in e} b_{e,i,j}( (v_s+w_s-n\delta_{s,l+1})_{s \in e \backslash \{i\}}, \sum_{s \in e} v_s + \sum_{s \in e} w_s - n + m_j )$$
for all $v = (v_1,\ldots,v_{l+1}) \in \Z_P^l$, $w = (w_1,\ldots,w_l)$, and $n \in [N] \cup [N']$, where we adopt the conventions that $v_{l+1} := -\Sigma(v_l)$, $w_{l+1} := -\Sigma(w)$, and $\delta_{s,l+1}$ is the Kronecker delta, equal to $1$ when $s=l+1$ and $0$ otherwise.

Now since $N, N' \leq F(M)$ and $M < M^*$ we see that $N, N' < M_k^{1/4}$.   Thus we see that $\sum_{s \in e} w_s - n = O( M_k^{1/2} )$.  On the other hand, $M_j \geq M_k$.  Thus we can shift $m$ by $\sum_{s \in e} w_s - n$ and write
$$ \varphi_{e,j}(v+w,-\Sigma(v+w)-n) = \E_{m_j \in [M_j]} \psi_{e,j,v,m_j}(w,-\Sigma(w)-n) + O( M_k^{-1/2} )$$
where $\psi_{e,j,v,m_j}: \Z^{l+1} \to [-1,1]$ is the function
$$ \psi_{e,j,v,m_j}(u) := \prod_{i \in e} b_{e,i,j}( (v_s+u_s)_{s \in e \backslash \{i\}}, \sum_{s \in e} v_s + m_j ).$$
From Definition \ref{complex-def} we observe that $\psi_{e,j,v,m_j}$ is a basic function of complexity at most $d-1$.

Applying \eqref{hek} (and recalling that $\Psi_e$ has degree and coefficients $O_\eps(1)$), we can now write
$$ h(v+w,-\Sigma(v+w)-n) = \prod_{e \in \I} \Psi_e\left( 
\left(\E_{m_j \in [M_j]} \psi_{e,j,v,m_j}(w,-\Sigma(w)-n)\right)_{j=k}^K \right) + O_\eps( M_k^{-1/2} ).$$
We expand out the polynomials $\Psi_e$ and collect all the $m_j$ averages, and eventually rewrite the right-hand side in the form
$$ \E_{\vec m \in \M} \sum_{\alpha \in A} c_\alpha \Phi_{v,\vec m, \alpha}(w, -\Sigma(w)-n) + O_\eps( M_k^{-1/2} ),$$
where $\M$ is a finite index set (it is the product of finitely many intervals of the form $[M_j]$), $A$ is another index set of size $O_\eps(1)$, the coefficients $c_\alpha$ are numbers of size $O_\eps(1)$, and the $\Phi_{v,\vec m,\alpha}: \Z^{l+1} \to [-1,1]$ are various basic functions of complexity at most $d-1$ whose exact form is not of importance to us (they are products of various $\psi_{e,j,v,m_j}$, where the $m_j$ are drawn from components of the $\vec m$).  If we then define
$$ \varphi_{v,\vec m} := \sum_{\alpha \in A} c_\alpha \Phi_{v,\vec m, \alpha}(w, -\Sigma(w)-n)$$
we obtain the claim.
\end{proof}

From this lemma and \eqref{evil}, we can now bound the left-hand side of \eqref{hon} by
$$
\left( \E_{v \in \Z_P^l, w \in [L]^l} \left|\E_{\vec m \in \M} \left(\E_{n \in [N]} \varphi_{v,\vec m}(w, -\Sigma(w)-n) - \E_{n \in [N']} \varphi_{v,\vec m}(w, -\Sigma(w)-n)\right)\right|^2 \right)^{1/2} + O_\eps( M_k^{-1/4} )$$
which by Cauchy-Schwarz can be bounded by
$$
\ll \left( \E_{v \in \Z_P^l, \vec m \in \M, w \in [L]^l} \left|\E_{n \in [N]} \varphi_{v,\vec m}(w, -\Sigma(w)-n) - \E_{n \in [N']} \varphi_{v,\vec m}(w, -\Sigma(w)-n)\right|^2 \right)^{1/2} + O_\eps( M_k^{-1/4} ).$$
The next step is to move from $[L]^l$ to a cyclic group.  Let $Y$ be the finite set $\Z_P^l \times \M$, which we endow with the uniform measure $\mu_Y$.  Let $Q := (l+1) L$.  
We define the functions $\tilde \varphi: \Z_Q^{l+1} \times Y \to [-1,1]$ by defining
$$ \tilde \varphi((w_1,\ldots,w_l, w_{l+1}), (v,\vec m)) := \varphi_{v,\vec m}(w_1,\ldots,w_l,w_{l+1})$$
when $v \in \Z_P^l$, $\vec m \in \M$, $w_1,\ldots,w_l \in [L]$ and $w_{l+1} \in \{-1,\ldots,-Q\}$ (where we identify integers with elements of $\Z_Q$ in the usual manner), and $\tilde \varphi_{v,\vec m}=0$ otherwise.  Note that as $\varphi_{v,\vec m}: \Z^{l+1} \to \R$ is an elementary function of complexity at most $(d-1,O_\eps(1))$, the function $\tilde \varphi: \Z_Q^{l+1} \times X \to \R$ is also.

Since $|Q^{l}| \ll |[L]^l|$, one can bound the left-hand side of \eqref{hon} by
$$ \ll \left( \E_{y \in Y, w \in \Z_Q^l} \left|\E_{n \in [N]} \tilde \varphi((w, -\Sigma(w)-n), y) - \E_{n \in [N']} \tilde \varphi((w, -\Sigma(w)-n), y)\right|^2 \right)^{1/2} + O_\eps( M_k^{-1/4} ),$$
which by Definition \ref{dap} can be expressed as
$$ \ll \| \Delta_N \tilde \varphi - \Delta_{N'} \tilde \varphi \|_{L^2( \Z_Q^{l+1} \times Y )} + O_\eps( M_k^{-1/4} ),$$
where we have abused notation slightly and viewed $\Delta_N$ as a function on $\Z_Q^{l+1} \times Y$ instead of $\Z_Q^l \times Y$ by adding a dummy variable.
But we can now apply the inductive hypothesis, Theorem \ref{main3}, with $d$, $P$, $X$, $g$, $\eps$, $M_*$, $J$ replaced by $d-1$, $Q$, $Y$, $\varphi$, $\eps/C$, $M_{**}$, and $O_\eps(1)$ respectively for some large absolute constant $C$, and conclude the existence of $M_{**} \leq M \ll O_{F,\eps,C,M_{**}}(1)$ such that
$$ \| \Delta_N \tilde \varphi - \Delta_{N'} \tilde \varphi \|_{L^2( \Z_Q^{l+1} \times Y )} \leq \eps/C$$
for all $M \leq N, N' \leq F(M)$.  If we choose $\tilde F$ to be sufficiently fast-growing depending on $F$, $\eps$, $C$, we can ensure that $M^{**} \geq M$.  The left-hand side of \eqref{hon} is now bounded by
$$ \ll \eps/C + O_\eps( M_k^{-1/4} ).$$
By making $C$ sufficiently large, and making $\tilde F$ sufficiently fast-growing depending on $\eps$, we thus establish \eqref{hon}.  This establishes Theorem \ref{main3b}, and hence (by induction) Theorem \ref{main3}.  Theorem \ref{main2} and Theorem \ref{main} then follow.
\hfill$\square$

\appendix

\section{A quantitative dominated convergence theorem}

We recall the following version case of the Lebesgue dominated convergence theorem on the net $\N^2$:

\begin{theorem}[Lebesgue dominated convergence theorem for $\N^2$]\label{ldct-g}  Let $(X,\X,\mu)$ be a probability space, and for each $n,n' \in \N$ let $f_{n,n'}: X \to [0,1]$ be a measurable function.  If, for almost every $x \in X$, we have $\lim_{n,n' \to \infty} f_{n,n'}(x) = 0$ then we have $\lim_{n, n' \to \infty} \int_X f_{n,n'}(x)\ d\mu(x) = 0$.
\end{theorem}

In this appendix we apply a correspondence principle (essentially the Furstenberg correspondence principle) to transfer this infinitary theorem to a finitary counterpart, which may be of some independent interest.  More precisely, we have

\begin{theorem}[Finitary Lebesgue dominated convergence theorem]\label{fin-ldct}  Suppose we have a positive integer $M_{*,F,\eps}$ assigned to each $\eps > 0$ and each function $F: \N \to \N$.  Then for every $\eps' > 0$ and every $F': \N \to \N$ we can find a positive integer $M'_{*,F',\eps'}$ with the following property: given any probability space $(X,\X,\mu)$, and sequence $f_{n,n'}: X \to [0,1]$ of measurable functions with the quantitative convergence property
\begin{itemize}
\item[(*)] For every $\eps > 0$ and every $F: \N \to \N$, for almost every $x \in X$ there exists an integer $1 \leq M \leq M_{*,F,\eps}$ such that $f_{n,n'}(x) \leq \eps$ for all $M \leq n,n' \leq F(M)$.
\end{itemize}
there exists an integer $1 \leq M \leq M'_{*,F',\eps'}$ such that
$$ \int_X f_{n,n'}\ d\mu \leq \eps'$$
for all $M \leq n,n' \leq F'(M)$.
\end{theorem}

\begin{remark} In this theorem, the indices $n,n'$ are ranging over the net $\N^2$, but it will be clear from the proof that one could in fact work with any countable net.
\end{remark}

\begin{proof}  Let us fix the assignment $(\eps,F) \mapsto M_{*,F,\eps}$, as well as the quantity $\eps' > 0$ and the function $F' > 0$.  We may assume that $F'(M) \geq M$ for all $M$ since the claim is vacuous otherwise.  Suppose for contradiction that the theorem failed for these parameters.  Untangling all the quantifiers carefully (and using the axiom of choice), this means that for every integer $m$ we can find a probability space $(X^{(m)}, \X^{(m)}, \mu^{(m)})$ and a family $f^{(m)}_{n,n'}: X^{(m)} \to [0,1]$ of sequences obeying the property (*), but such that
\begin{equation}\label{mnnm}
 \sup_{M \leq n, n' \leq F'(M)} \int_{X^{(m)}} f^{(m)}_{n,n'}\ d\mu^{(m)} > \eps
\end{equation}
for all $1 \leq M \leq m$.

Let $[0,1]^{\N^2}$ be the space of all functions $g: \N \times \N \to [0,1]$; by Tychonoff's theorem, this is a compact Hausdorff topological space with the product topology and the usual Borel $\sigma$-algebra $\B$, which is countably generated. 

Observe that we have the maps $f^{(m)}: X^{(m)} \to [0,1]^{\N^2}$ for each $m \geq 1$ defined by 
$$ f^{(m)}(x)(n,n') := f^{(m)}_{n,n'}(x).$$
One easily verifies that this map is measurable.  Thus, we can push forward the probability measure $\mu^{(m)}$ forward by $f^{(m)}$ to create a probability measure $\nu^{(m)} := f^{(m)}_* \mu^{(m)}$ on $[0,1]^{\N^2}$. The space of probability measures on the countably generated $\sigma$-algebra $\B$ is weakly sequentially compact.  What this means is that we can find a subsequence $\nu^{(m_j)}$ of the probability measures $\nu^{(m)}$ which converge weakly to another probability measure $\nu$ on $[0,1]^{\N^2}$ in the sense that 
\begin{equation}\label{numj}
\lim_{j \to \infty} \nu^{(m_j)}(A) = \nu(A)
\end{equation}
for any elementary set $A$.  Indeed, for each elementary set $A$ one can refine the subsequence $m_j$ so that $\nu^{(m_j)}(A)$ is convergent, and then by the usual Arzela-Ascoli type diagonalisation argument we can ensure that $\nu^{(m_j)}(A)$ converges to a limit $\nu(A)$ for all elementary sets $A$.  One can then use the Caratheodory extension theorem or Kolmogorov extension theorem to extend $\nu$ to a probability measure.

Fix this subsequence $m_j$ and the limit measure $\nu$.  For any natural numbers $n,n' \in \N$, let $\pi_{n,n'}: [0,1]^{\N^2} \to [0,1]$ be the coordinate projection $\pi_{n,n'}(g) := g(n,n')$.  These functions are continuous on $[0,1]^{\N^2}$ and hence measurable; indeed we see that $\pi_{n,n'}^{-1}([a,b])$ is an elementary set for any interval $[a,b]$ with rational endpoints.  From
\eqref{mnnm} and the definition of $\nu^{(m)}$ and $\pi_{n,n'}$ we see that
$$
 \sup_{M \leq n, n' \leq F'(M)} \int_{[0,1]^{\N^2}} \pi_{n,n'}(y)\ d\nu^{(m)}(y) > \eps'
$$
for all $1 \leq M \leq m$.  Fixing $M$, specialising $m$ to $m_j$ for $j$ sufficiently large, and then taking limits as $j \to \infty$ using the weak convergence of the $\nu^{(m_j)}$ (noting that the level sets of $\pi_n-\pi_{n'}$ are elementary sets), we conclude that
$$
 \sup_{M \leq n, n' \leq F'(M)} \int_{[0,1]^{\N^2}} \pi_{n,n'}(y)\ d\nu(y) \geq \eps'
$$
for all $M \geq 1$.  Observe that the function $(n,n') \mapsto \pi_{n,n'}(y)$ is (tautologically) a pseudometric taking values in $[0,1]$ for each $y \in [0,1]^{\N^2}$.  Applying Theorem \ref{ldct-g} in the contrapositive, we conclude that
$$ \nu( \{ y \in [0,1]^{\N^2}: \inf_{M \to \infty} \sup_{n,n' \geq M} \pi_{n,n'}(y) > 0 \} ) > 0.$$
By countable subadditivity\footnote{This can be viewed as an infinite version of the pigeonhole principle, viz. if a set of positive measure is covered by countably many measurable sets, then at least one of those sets also has positive measure.}, this implies that there exists an $\eps > 0$ such that
$$ \nu( \{ y \in [0,1]^{\N^2}: \sup_{n,n' \geq M} \pi_{n,n'}(y) \geq 2\eps \hbox{ for all } M \geq 1 \} ) > 0.$$
If we define the sets
$$
E_{M,M'} := \{ y \in [0,1]^{\N^2}: \sup_{M \leq n,n' \leq M'} \pi_{n,n'}(y) \geq 2\eps \}
$$
for all $1 \leq M \leq M'$, we thus see that
$$ \nu( \bigcap_{M=1}^\infty \bigcup_{M'=M}^\infty E_{M,M'} ) > 0.$$
By using countable subadditivity recursively, we can thus find an integer $F(M) \geq M$ associated to every $M \geq 1$ such that
$$ \nu( \bigcap_{M=1}^{M_0} E_{M,F(M)} \cap \bigcup_{M=M_0+1}^\infty \bigcup_{M'=M}^\infty E_{M,M'} ) > 0$$
for all $M_0 \geq 1$, and in particular that
\begin{equation}\label{emfm}
 \nu( \{ y \in [0,1]^{\N^2}: \inf_{1 \leq M \leq M_0} \sup_{M \leq n,n' \leq F(M)} \pi_{n,n'}(y) \geq 2\eps \} ) > 0
\end{equation}
for all $M_0 \geq 1$.

Fix this $F$.  We apply hypothesis (*) for the sequences $f^{(m_j)}_{n,n'}(x)$ and conclude that for all $j$ and $\mu^{(m_j)}$-almost every $x \in X^{(m_j)}$ we have
$$ \inf_{1 \leq M \leq M_{*,F,\eps}} \sup_{M \leq n,n' \leq F(M)} f^{(m_j)}_{n,n'}(x) \leq \eps.$$
Equivalently, by the definition of $\nu^{(m_j)}$ and $\pi_n$ we have
$$ \nu^{(m_j)}(\{ y \in [0,1]^{\N^2}: \inf_{1 \leq M \leq M_{*,F,\eps}} \sup_{M \leq n,n' \leq F(M)} \pi_{n,n'}(y) \leq \eps \}) = 1.$$
Since $\nu^{(m_j)}$ converges weakly to $\nu$, and the subset of $[0,1]^{\N^2}$ appearing above is compact and depends on only finitely many coordinates of $[0,1]^{\N^2}$, we conclude that
$$ \nu(\{ y \in [0,1]^{\N^2}: \inf_{1 \leq M \leq M_{*,F,\eps}} \sup_{M \leq n,n' \leq F(M)} \pi_{n,n'}(y) \leq \eps \}) = 1.$$
But this contradicts \eqref{emfm}.  The proof of Theorem \ref{fin-ldct} is complete.
\end{proof}

\begin{remark} In principle, the quantity $M'_{*,F',\eps'}$ can be explicitly computed from $F'$, $\eps'$, and the assigment $(F,\eps) \mapsto M_{*,F,\eps}$.  In practice, though, it seems remarkably hard to do; the proof of the Lebesgue dominated convergence theorem, if inspected carefully, relies implicitly on the infinite pigeonhole principle, which is notoriously hard to finitise.  Indeed the situation here is somewhat reminiscent of that of the Paris-Harrington theorem \cite{paris}.  Note that it was established in \cite{yu} (see also \cite{simic}) that the Lebesgue dominated convergence theorem is equivalent in the reverse mathematics sense to the arithmetic comprehension axiom $(ACA)$, which does strongly suggest that the dependence of $M'_{*,F',\eps'}$ on the above parameters is likely to be fantastically poor.
\end{remark}

We will use Theorem \ref{fin-ldct} to eliminate the role of various probability spaces $(X,\X,\mu)$ in our analysis.  This elimination is not, strictly speaking, absolutely necessary\footnote{This is analogous to how the ergodic decomposition is, strictly speaking, not necessary in the proof of many recurrence theorems in ergodic theory, if one is willing to make all of one's computations ``relative'' to the shift-invariant factor.} for us; we could instead passively carry such spaces with us throughout our arguments, at the cost of making the notation in those arguments slightly more complicated.  We have however chosen this approach to highlight the finitary version of the dominated convergence theorem, which is not so well-known in the literature.

\end{document}